\documentclass[11pt,a4paper]{amsart}
\usepackage[english]{babel}
\usepackage[T1]{fontenc}
\usepackage{mathpazo,amssymb}
\usepackage{upref}
\usepackage{enumerate}
\usepackage{graphicx,xcolor}
\usepackage{cancel}
\usepackage{dsfont}       
\usepackage[colorlinks=true, pdftitle={Regularizing effect homogeneous
  evolution equations}, pdfauthor={Daniel Hauer
and Jos\'e M. Maz\'on} ]{hyperref}

\title[Homogeneous evolution equations]{Regularizing effect of homogeneous evolution equations : Case homogeneous order
zero}

\author{Daniel Hauer} \address[Daniel Hauer]{School of Mathematics and
  Statistics, The University of Sydney, NSW 2006, Australia}
\email{\href{mailto:daniel.hauer@sydney.edu.au}{\nolinkurl{daniel.hauer@sydney.edu.au}}}

\author{Jos\'e M. Maz\'on} \address[Jos\'e M. Maz\'on]{Departament
  d'An\`alisi Matem\`atica, Universitat de Val\`encia, Valencia,
  Spain} \email{\href{mailto:mazon@uv.es}{\nolinkurl{mazon@uv.es}}}
\thanks{The second author have been partially supported by the Spanish
  MINECO and FEDER, project MTM2015-70227-P. The first author is very
  grateful for the kind invitation to the Universitat de Val\`encia
  and their hospitality.}



\numberwithin{equation}{section}

\newtheorem{theorem}{Theorem}[section]
\newtheorem{proposition}[theorem]{Proposition}
\newtheorem{lemma}[theorem]{Lemma}
\newtheorem{corollary}[theorem]{Corollary}

\theoremstyle{definition}
\newtheorem{definition}[theorem]{Definition}

\newtheorem{remark}[theorem]{Remark}
\newtheorem{example}[theorem]{Example}

\newcommand\R{{\mathbb{R}}}
\newcommand\N{\mathbb{N}}

\newcommand\W{\mathcal{W}}

\newcommand\dx{\mathrm{d}x }
\newcommand\dy{\mathrm{d}y }
\newcommand\dr{\mathrm{d}r }
\newcommand\ds{\mathrm{d}s }

\newcommand\dmu{\mathrm{d}\mu}
\newcommand\dt{\mathrm{d}t }
\newcommand\td{\mathrm{d} }

\newcommand\abs[1]{\lvert#1\rvert}
\newcommand\labs[1]{\left\lvert#1\right\rvert}
\newcommand\norm[1]{\lVert#1\rVert}
\newcommand\lnorm[1]{\left\lVert#1\right\rVert}

\definecolor{darkred}{rgb}{0.7,0.1,0.1}

\begin{document}
\date{\today}
\maketitle

\tableofcontents


\begin{abstract}
  In this paper, we develop a functional analytical theory for establishing that mild
  solutions of first-order Cauchy problems involving homogeneous
  operators of order zero are strong solutions; in particular, the first-order
  time derivative satisfies a global regularity estimate depending
  only on the initial value and the positive time. We apply those
  results to the Cauchy problem associated with the total variational
  flow operator and the nonlocal fractional $1$-Laplace operator.
\end{abstract}

%
%

\section{Introduction}

In the pioneering work~\cite{MR648452}, B\'enilan and Crandall showed that
for the class of \emph{homogeneous operators $A$ of order $\alpha>0$
with $\alpha\neq 1$}, defined on a normed space
$(X,\norm{\cdot}_{X})$, every solution of the differential inclusion
\begin{equation}
  \label{eq:10}
  \frac{\td u}{\dt}+A(u(t))\ni 0
\end{equation}
satisfies the global regularity estimate
\begin{equation}
  \label{eq:13}
  \limsup_{h\to 0+}\frac{\norm{u(t+h)-u(t)}_{X}}{h}\le
  2\,L\frac{\norm{u_{0}}_{X}}{\abs{\alpha-1}}\frac{1}{t}\qquad\text{for
    every $t>0$.}
\end{equation}
Here, $A\subseteq X\times X$ might be multi-valued and is called
\emph{homogeneous of order $\alpha$} if
\begin{equation}
  \label{eq:41}
  A(\lambda u)=\lambda^{\alpha}Au\qquad\text{ for all $\lambda\ge 0$
and $u\in D(A)$.}
\end{equation}
Moreover, to obtain~\eqref{eq:13}, it is assumed that there is a family
$\{T_{t}\}_{t\ge 0}$ associated with $A$ of Lipschitz continuous
mappings $T_{t}$ on $X$ of constant $L$ such that  
\begin{equation}
   \label{eq:22}
   u(t)=T_{t}u_{0}\qquad\text{for every $t\ge 0$,}
\end{equation}
is (in some given sense) a solution of~\eqref{eq:10} for some initial
value $u_{0}\in X$. We refer to Definition~\ref{def:mild-solution} and
Definition~\ref{def:strong-sols} for the different notions of solutions.

Further, if $X$ is equipped with a partial ordering
$``\!\!\!\le''$ such that $(X,\le)$ defines an ordered vector space,
and if for this ordering, the family $\{T_{t}\}_{t\ge 0}$ is
\emph{order-preserving} (that is,~\eqref{eq:23} below holds), then
every \emph{positive}\footnote{Here, we call a measurable function $u$
  \emph{positive} if $u\ge 0$ for the given partial ordering
  $``\!\!\!\le''$.} solution $u$ of~\eqref{eq:10} satisfies
the 
point-wise estimate
\begin{equation}
  \label{eq:14}
  (\alpha-1)\frac{\td u}{\dt}_{\!\!+}\!\!(t)\ge - \frac{u}{t}
  \qquad\text{in $\mathcal{D}'$ for every $t>0$.}
\end{equation}

Estimates of the form~\eqref{eq:13} describe an instantaneous and
global \emph{regularizing effect} of solutions $u$ of~\eqref{eq:10},
since they imply that the solution $u$ of~\eqref{eq:10} is locally
Lipschitz continuous in $t\in(0,+\infty)$. Further~\eqref{eq:14}
provides a rate of dissipativity involved in the differential
inclusion~\eqref{eq:10}.\medskip

It is the aim of this paper to extend the theory developed
in~\cite{MR648452} to the important case $\alpha=0$; in other words,
for the class of \emph{homogeneous operators $A$ of order zero} (see
Definition~\ref{def:1} below). Important examples of this class of
operators include the (negative) \emph{total variational flow operator}
$Au=-\Delta_{1}u:=-\textrm{div}\left(\frac{D u}{\abs{D u}}\right)$, also known as
(negative) \emph{$1$-Laplacian}, or the \emph{$1$-fractional Laplacian}
\begin{displaymath}
	Au=(-\Delta_{1})^{s}u(x):=PV\int_{\Sigma}\frac{u(y)-u(x)}{\abs{u(y)-u(x)}}
        \frac{\dy}{\abs{x-y}^{d+s}}\;\mbox{},
	\qquad s\in (0,1).
\end{displaymath}
In our first main result (Theorem~\ref{thm:1bis}), we establish the
global regularity estimate~\eqref{eq:13} for order
$\alpha=0$ and for solutions $u$ of \emph{differential inclusions
  with a forcing term}:
\begin{equation}
  \label{eq:10bis}
  \frac{\td u}{\dt}+A(u(t))\ni f(t)\qquad\text{on $(0,T)$,}
\end{equation}
where $f : [0,T]\to X$ is an integrable function, and $T>0$. In
Corollary~\ref{thm:1} and Corollary~\ref{cor:1}, we provide the
resulting inequality when $f\equiv 0$ and the right hand-side
derivative $\frac{\td u}{\dt_{\!+}}(t)$ of $u$ exists at $t>0$.

In many applications (cf~Section~\ref{sec:application}), $X$ is given
by the classical Lebesgue space $(L^{r},\norm{\cdot}_{r})$,
($1\le r\le \infty$). If $\{T_{t}\}_{t\ge 0}$ 
is a semigroup satisfying an \emph{$L^{q}$-$L^{r}$-regularity estimate}
\begin{equation}
  \label{eq:44}
  \norm{T_{t}u_{0}}_{r}\le
  C\,e^{\omega t}\frac{\norm{u_{0}}_{q}^{\gamma}}{t^{\delta}}\qquad
  \text{for all $t>0$, and $u_{0}\in L^{q}$,}
\end{equation}
for $\omega\in \R$,
$\gamma=\gamma(q,r,d)$, $\delta=\delta(q,r,d)>0$, and some (or
for all) $1\le q<r$, then we show in Corollary~\ref{cor:1} that
combining~\eqref{eq:13} with~\eqref{eq:44} yields
\begin{equation}
  \label{eq:16}
  \limsup_{h\to 0+}\frac{\norm{u(t+h)-u(t)}_{r}}{h}
  \le C\,L\,2^{\delta+2}\,e^{\omega\,t}
  \frac{\norm{u_{0}}_{q}^{\gamma}}{t^{\delta+1}}.
\end{equation}
Regularity estimates similar to~\eqref{eq:44} have been studied
recently by many authors (see, for
example,~\cite{MR1103113,MR1218884,MR2569498} covering the linear
theory, and~\cite{CoulHau2017} the nonlinear one and there references
therein).

In Theorem~\ref{thm:2}, Corollary~\ref{cor:2}
and Corollary~\ref{cor:3} we generalize the point-wise
estimate~\eqref{eq:14} to the homogenous order $\alpha=0$.

We emphasize that the regularizing effect of solutions $u$
of~\eqref{eq:10} remains true with a slightly different inequality
(see Corollary~\ref{cor:1bisF}) if the homogeneous operator $A$ is
perturbed by a Lipschitz mapping $F$. This is quite surprising since
$F$ might not be homogenous and hence, the operator $A+F$ is also not
homogeneous.

In Section~\ref{sec:semigroups}, we consider the class of
\emph{quasi accretive operators} $A$ (see
Definition~\ref{def:quasi-accretive}) and outline how the property
that $A$ is homogeneous of order zero is passed on to the
\emph{semigroup} $\{T_{t}\}_{t\ge 0}$ generated by $-A$ (see the
paragraph after Definition~\ref{def:mild-solution}). In particular, we
discuss when solutions $u$ of~\eqref{eq:10} are differentiable a.e. in
$t>0$.

The fact that every Lipschitz continuous mappings $u : [0,T]\to X$ is
differentiable almost everywhere on $(0,T)$ depends on the underlying
geometry of the given Banach space $X$; this property is well-known as
the Radon-Nikod\'ym property of a Banach space. The Lebesgue space
$L^{1}$ has not this property, but alone from the physical point of
view, $L^1$ is for many models not avoidable. In~\cite{MR1164641},
B\'enilan and Crandall developed the celebrated theory of
\emph{completely accretive} operators $A$ (in $L^{1}$). For this class
of operators, it is known that for each solutions $u$ of~\eqref{eq:10}
in $L^1$, the derivative $\frac{\td u}{\dt}$ exists in $L^{1}$. These
results have been extended recently to the notion of \emph{quasi
  completely accretive} operators in~\cite{CoulHau2017}. In
Section~\ref{sec:completely-accretive}, we study regularity
estimates of the form~\eqref{eq:13} for $\alpha=0$ satisfied by solutions $u$
of~\eqref{eq:10}, where $A$ is a quasi completely accretive operator
of homogeneous order zero. In fact, the two operators $-\Delta_{1}$ and
$(-\Delta_{1})^{s}$ mentioned above, belong exactly to this class of
operators. Thus, our two main examples of differential inclusions
discussed in Section~\ref{sec:application} are
\begin{align}\label{eq:15}
  & \frac{\td u}{\dt}-\textrm{div}\left(\frac{D u}{\abs{D
    u}}\right)+f(\cdot,u)\ni 0\hspace{1.6cm}
    \\[7pt]
\label{eq:17}
  &  \frac{\td
    u}{\dt}+PV\int_{\Sigma}\frac{(u(y)-u(x))}{\abs{u(y)-u(x)}} \frac{\dy}{\abs{x-y}^{d+s}}+f(\cdot,u)\ni
    0, 
\end{align}
respectively equipped with some boundary conditions on a domain
$\Sigma$ in $\R^{d}$, $d\ge 1$. In~\eqref{eq:15} and~\eqref{eq:17},
the function $f$ is a Carath\'eo\-dory function,
which is Lipschitz continuous in the second variable with constant $\omega>0$
uniformly with respect to the first variable (see
Section~\ref{sec:application} for more details).\medskip

Note , if the right hand-side
derivative $\frac{\td u}{\dt_{\!+}}(t)$ of a solution $u$ of~\eqref{eq:10} exists at every $t\in
(0,1]$, then~\eqref{eq:13} for $\alpha=0$ becomes
\begin{equation}
  \label{eq:49}
  \norm{Au(t)}_{X}\le 2\,L\frac{\norm{u_{0}}_{X}}{t}\qquad\text{for
    every $t>0$.}
\end{equation}
Here, it is worth mentioning that if the operator $A$ in~\eqref{eq:10}
is linear (that is, $\alpha=1$), then inequality~\eqref{eq:49} means
that $-A$ generates an analytic semigroup $\{T_{t}\}_{t\ge 0}$
(cf~\cite{MR2103696,MR710486}). Thus, it is interesting to see that a
similar regularity inequality such as~\eqref{eq:49}, in particular, holds for certain
classes of nonlinear operators. In addition, if
$\norm{\cdot}_{X}$ is the induced norm by an inner product $(\cdot,\cdot)_{X}$ of a
Hilbert space $X$ and $A$ is a sub-differential operators
$\partial\varphi$ on $X$, then inequality~\eqref{eq:49} is also
satisfied by solutions of~\eqref{eq:10} (cf~\cite{MR0348562}). In
\cite{MR1799898}, inequality~\eqref{eq:49} was shown to hold for
solutions of~\eqref{eq:15} with $f\equiv 0$ and equipped with Neumann
boundary conditions. 


%
%

\section{Main results}
\label{sec:main-results}

Suppose $X$ is a linear vector space and $\norm{\cdot}_{X}$ a
semi-norm on $X$. Then, the main object of this paper is the following
class of operators.

\begin{definition}\label{def:1}
  An operator $A$ on $X$ is said to be \emph{homogeneous
    of order zero} if for every $u\in D(A)$ and $\lambda\ge 0$, one has
  that $\lambda u\in D(A)$, and $A$ satisfies~\eqref{eq:41} for $\alpha=0$.
  \end{definition}

  \begin{remark}
    It follows necessarily from~\eqref{eq:41} that for every
    homogeneous operator $A$ of order $\alpha>0$, one has that
    $0\in A0$. But for homogeneous operators $A$ of order zero, the
    property $0\in A0$ does not need to hold.
  \end{remark}

  Now, assume that for the operator $A$ on $X$ and for given
  $f : [0,T]\to X$ and $u_{0}\in X$, the function
  $u\in C^{1}([0,T];X)$ is a classical solution of the differential
  inclusion~\eqref{eq:10bis} with forcing term $f$ satisfying
  \emph{initial value} $u(0)=u_{0}$. If $A$ is homogeneous of order
  zero, then for $\lambda>0$, the function
\begin{displaymath}
  v(t)=\lambda^{-1} u(\lambda t),\qquad( t\in [0,T]),
\end{displaymath}
satisfies
\begin{displaymath}
  \frac{\td v}{\dt}(t)=\frac{\td u}{\dt}(\lambda t)\in
  - A(u(\lambda t))+f(\lambda t)
  =-A(v(t)) +f(\lambda t)
\end{displaymath}
for every $t\in (0,T)$ with initial value
$v(0)=\lambda^{-1} u(0)=\lambda^{-1} u_{0}$. Thus, if for every $t\in
[0,T]$, we denote
\begin{equation}
  \label{eq:36}
  T_{t}(u_{0},f):=u(t)\qquad\text{for every $u_{0}$ and $f$,}
\end{equation}
where $u$ is the unique classical solution $u$ of~\eqref{eq:10bis} with
initial value $u(0)=u_{0}$, then the above reasoning shows that the homogeneity of $A$
is reflected in
\begin{equation}
 \label{eq:21bis}
  \lambda^{-1}T_{\lambda
    t}(u_{0},f)=T_{t}(\lambda^{-1}u_{0},f(\lambda\cdot))\qquad
  \text{for every $\lambda>0$,}
\end{equation}
and all $t\in [0,T]$. Identity~\eqref{eq:21bis} together with standard
growth estimates of the form
\begin{equation}
  \label{eq:20bis}
  \begin{split}
    &e^{-\omega t}\norm{T_{t}(u_{0},f)-T_{t}(\hat{u}_{0},\hat{f})}_{X}\\
    &\qquad\le L e^{-\omega
      s}\norm{T_{s}(u_{0},f)-T_{s}(\hat{u}_{0},\hat{f})}_{X}
    +L \int_{s}^{t}e^{-\omega
      r}\,\norm{f(r)-\hat{f}(r)}_{X}\,\dr
  \end{split}
\end{equation}
for every $0 \le s \le t (\le T)$, (for some $\omega\in \R$ and
$L\ge 1$) are the main ingredients to obtain global regularity
estimates of the form~\eqref{eq:13}. This leads to our first main
result.

\begin{theorem}\label{thm:1bis}
  For a subset $C\subseteq X$, let $\{T_{t}\}_{t=0}^{T}$ be a family
  of mappings $T_{t} : C\times L^{1}(0,T;X)\to C$
  satisfying~\eqref{eq:20bis}, \eqref{eq:21bis}, and
  $T_{t}(0,0)\equiv 0$ for all $t\ge 0$. Then for every $u_{0}\in C$,
  $f\in L^{1}(0,T;X)$, and $t\in (0,T]$, $h>0$, one has that
\begin{equation}
  \label{eq:25bis}
  \begin{split}
    &\norm{T_{t+h}(u_{0},f)-T_{t}(u_{0},f)}_{X}\\
    &\qquad \le \tfrac{\abs{h}}{t}\;L\,e^{\omega\,t}\left[2
  \norm{u_{0}}_{X}+ \left(1+\tfrac{h}{t}\right) \,\int_{0}^{t}e^{-\omega
    s}\lnorm{\frac{f(s+\tfrac{h}{t}s)-f(s)}{\frac{h}{t}}}_{X}\,\ds\right.\\
  &\hspace{8cm}\left.+\int_{0}^{t}e^{-\omega
    s}\norm{f(s)}_{X}\,\ds\right].
  \end{split}
\end{equation}
In particular, if
\begin{displaymath}
  V(f,t):=\limsup_{\xi\to 0} \int_{0}^{t}e^{-\omega s}\lnorm{\frac{f(s+\xi s)-f(s)}{\xi}}_{X}\,\ds,
\end{displaymath}
then the family $\{T_{t}\}_{t\ge 0}$ satisfies
\begin{equation}
  \label{eq:38}
  \begin{split}
    &\limsup_{h\to 0+}\lnorm{\frac{T_{t+h}(u_{0},f)-T_{t}(u_{0},f)}{h}}_{X}\\
    & \qquad \le \frac{L e^{\omega\,t}}{t} \left[2
  \norm{u_{0}}_{X}+V(f,t) +\int_{0}^{t}e^{-\omega
    s}\norm{f(s)}_{X}\,\ds\right].
  \end{split}
\end{equation}
for every $t>0$, $u_{0}\in C$, $f\in L^{1}(0,T;X)$, and if $f$ is
locally absolutely continuous and differentiable a.e. on $(0,T)$, then
\begin{equation}
  \label{eq:39}
  \begin{split}
    &\limsup_{h\to 0+}\lnorm{\frac{T_{t+h}(u_{0},f)-T_{t}(u_{0},f)}{h}}_{X}\\
    & \qquad \le \frac{L e^{\omega\,t}}{t}\, \left[2
  \norm{u_{0}}_{X}+ \int_{0}^{t}e^{-\omega
    s}s\norm{f'(s)}_{X}\,\ds+\int_{0}^{t}e^{-\omega
    s}\norm{f(s)}_{X}\,\ds\right].
  \end{split}
\end{equation}
Moreover, if 
  the right hand-side derivative
  $\frac{\td }{\dt_{+}}T_{t}(u_{0},f)$ exists (in $X$) at $t>0$, then
  \begin{equation}
    \label{eq:40}
      \lnorm{\frac{\td T_{t}(u_{0},f)}{\dt}_{\!\!+}}_{X} \le \frac{L e^{\omega\,t}}{t}\left[2
  \norm{u_{0}}_{X}+ V(t,f)+\int_{0}^{t}e^{-\omega
    s}\norm{f(s)}_{X}\ds\right].
  \end{equation}
\end{theorem}

\mbox{}\medskip

\allowdisplaybreaks
\begin{proof}
  Let $u_{0}\in C$, $f\in L^{1}(0,T;X)$, and for $t>0$, let $h\neq 0$
  satisfying $1+\frac{h}{t}\ge 0$. Then, choosing
  $\lambda=1+\frac{h}{t}$ in~\eqref{eq:21bis} gives
  \begin{equation}
    \label{eq:16bis}
    \begin{split}
      &T_{t+h}(u_{0},f)-T_{t}(u_{0},f)\\
      &\qquad =T_{\lambda
        t}(u_{0},f)-T_{t}(u_{0},f)\\
      &\qquad= \left(1+\tfrac{h}{t}\right)
      T_{t}\left[\left(1+\tfrac{h}{t}\right)^{-1}u_{0},f(\cdot+\tfrac{h}{t}\cdot)\right]-T_{t}(u_{0},f)
    \end{split}
  \end{equation}
  and so,
  \begin{equation}
  \begin{split}\label{eq:42}
    &T_{t+h}(u_{0},f)-T_{t}(u_{0},f)\\
    & \qquad =
        \left(1+\tfrac{h}{t}\right)\left[T_{t}\left[\left(1+\tfrac{h}{t}\right)^{-1}u_{0},f(\cdot+\tfrac{h}{t}\cdot)\right]
                                     -T_{t}(u_{0},f(\cdot+\frac{h}{t}\cdot))\right]\\
                     &\hspace{2cm} + \left(1+\tfrac{h}{t}\right)
                       \left[T_{t}\left[u_{0},f(\cdot+\tfrac{h}{t}\cdot)\right]-T_{t}(u_{0},f)\right]\\
                     &\hspace{4cm}    + \left[\left(1+\tfrac{h}{t}\right)-1\right]\, T_{t}(u_{0},f).
  \end{split}
\end{equation}
  Thus, by applying~\eqref{eq:20bis} and since $T_{t}(0,0)\equiv 0$,
  one sees that
\begin{align*}
  &\norm{T_{t+h}(u_{0},f)-T_{t}(u_{0},f)}_{X}\\
  &\qquad\le
  \left(1+\tfrac{h}{t}\right)
                                      \,\lnorm{T_{t}\left[\left(1+\tfrac{h}{t}\right)^{-1}u_{0},f(\cdot+\tfrac{h}{t}\cdot)\right]
                                     -T_{t}(u_{0},f(\cdot+\tfrac{h}{t}\cdot))}_{X}\\
  &\hspace{2cm} + \left(1+\tfrac{h}{t}\right)\lnorm{T_{t}(u_{0},f(\cdot+\tfrac{h}{t}\cdot))-T_{t}(u_{0},f)}_{X}\\
  &\hspace{4cm}    +
    \left[\left(1+\tfrac{h}{t}\right)-1\right]\,\norm{T_{t}(u_{0},f)}_{X}\\
  &\qquad\le\left(1+\tfrac{h}{t}\right)\;L\,e^{\omega\,t}\lnorm{\left(1+\tfrac{h}{t}\right)^{-1}u_{0}-u_{0}}_{X}\\
  &\hspace{2cm} +
    \left(1+\tfrac{h}{t}\right)\;L\,\int_{0}^{t}e^{\omega
    (t-s)}\norm{f(s+\tfrac{h}{t}s)-f(s)}_{X}\,\ds\\
  &\hspace{4cm}    + \;L\,e^{\omega\,t}
    \labs{\left(1+\tfrac{h}{t}\right)-1}\left(\norm{u_{0}}_{X}+\int_{0}^{t}e^{-\omega
    s}\norm{f(s)}_{X}\,\ds\right).
\end{align*}
From this is clear that~\eqref{eq:25bis}-\eqref{eq:40} follows.
\end{proof}

In the case $f\equiv 0$, then  the mapping $T_{t}$ given
by~\eqref{eq:36} only depends on the initial value $u_{0}$, that is,
\begin{equation}
  \label{eq:55}
  T_{t}u_{0}=T_{t}(u_{0},0)\qquad\text{for every $u_{0}$ and $t\ge 0$.}
\end{equation}
In this case, the estimates in Theorem~\ref{thm:1bis} reduce to the
following one.

\begin{corollary}\label{thm:1}
  Let $\{T_{t}\}_{t\ge 0}$ be a family of mappings $T_{t} : C\to C$
  defined on a subset $C\subseteq X$ satisfying
 \begin{align}
  \label{eq:20}
  \norm{T_{t}u_{0}-T_{t}\hat{u}_{0}}_{X}&\le L\,e^{\omega
    t}\,\norm{u _{0}-\hat{u}_{0}}_{X}\qquad\text{for all $t\ge 0$,
                                          $u$, $\hat{u}\in C$,}\\
   \label{eq:21}
  \lambda^{-1}\,T_{\lambda t}u _{0} &=T_{t}[\lambda^{-1}u _{0}]\qquad\text{for all
    $\lambda>0$, $t\ge 0$
    and $u _{0}\in C$,}
 \end{align}
 and $T_{t}0\equiv 0$ for all $t\ge 0$. Then, for every $u_{0}\in C$ and $t$, $h>0$, one has that
 \begin{equation}
   \label{eq:25}
   \norm{T_{t+h}u_{0}-T_{t}u_{0}}_{X}\le 2\,\tfrac{h}{t}\;L\,e^{\omega\,t}\norm{u_{0}}_{X}.
 \end{equation}
 In particular, the family $\{T_{t}\}_{t\ge 0}$ satisfies
 \begin{equation}
   \label{eq:18bis}
   \limsup_{h\to 0+}\frac{\norm{T_{t+h}u_{0}-T_{t}u_{0}}_{X}}{h}\le
   2 L e^{\omega t}\frac{\norm{u_{0}}_{X}}{t}\qquad\text{for every $t>0$, $u_{0}\in C$.}
 \end{equation}
 Moreover, if 
 the right hand-side derivative
 $\frac{\td }{\dt_{+}}T_{t}u_{0}$ exists (in $X$) at $t>0$, then
 \begin{equation}
   \label{eq:18bis2}
   \lnorm{\frac{\td T_{t}u_{0}}{\dt_{\!+}}}_{X}\le
   2 L e^{\omega t}\frac{\norm{u_{0}}_{X}}{t}.
 \end{equation}
\end{corollary}

For our next corollary, we recall the following well-known
definition.

\begin{definition}
  Let $C$ be a subset of $X$. Then, a family $\{T_{t}\}_{t\ge 0}$ of mappings $T_{t} : C\to C$ is
  called a \emph{semigroup} if
  \begin{math}
    T_{t+s}u=T_{t}\circ T_{s}u\text{ for every $t$, $s\ge0$, $u\in C$.}
  \end{math}
\end{definition}

\begin{corollary}\label{cor:1}
  Let $\{T_{t}\}_{t\ge 0}$ be a semigroup of mappings $T_{t} : C\to C$
  defined on a subset $C\subseteq X$ and suppose, there is a second
  vector space $Y$ with semi-norm $\norm{\cdot}_{Y}$ such that
  $\{T_{t}\}_{t\ge 0}$ satisfies the following $Y$-$X$-regularity estimate
  \begin{equation}\label{eq:9}
    \norm{T_{t}u_{0}}_{X}\le
    M\,e^{\hat{\omega} t}\frac{\norm{u_{0}}_{Y}^{\gamma}}{t^{\delta}}\qquad\text{for
      every $t>0$ and $u_{0}\in C$}
  \end{equation}
  for some $M$, $\gamma$, $\delta>0$ and
  $\hat{\omega}\in \R$. If for $u_{0}\in C$, $\{T_{t}\}_{t\ge 0}$
  satisfies~\eqref{eq:18bis}, then
  \begin{displaymath}
    \limsup_{h\to 0+}\frac{\norm{T_{t+h}u_{0}-T_{t}u_{0}}_{X}}{h}\le  2^{\delta+2} L\,
    M\,e^{\frac{1}{2}(\omega+\hat{\omega}) t}
      \frac{\norm{u_{0}}_{Y}^{\gamma}}{t^{\delta+1}}.
  \end{displaymath}
  Moreover, if the right hand-side derivative
 $\frac{\td }{\dt_{\!+}}T_{t}u_{0}$ exists (in $X$) at $t>0$, then
  \begin{displaymath}
    \lnorm{\frac{\td T_{t}u_{0}}{\dt _{+}}}_{X}\le  2^{\delta+2} L\,
    M\,e^{\frac{1}{2}(\omega+\hat{\omega}) t}
      \frac{\norm{u_{0}}_{Y}^{\gamma}}{t^{\delta+1}}.
  \end{displaymath}
\end{corollary}

\begin{proof}
  Since $\{T_{t}\}_{t\ge 0}$ is a semigroup, one sees
  by~\eqref{eq:18bis} and~\eqref{eq:9} that
  \begin{align*}
    \limsup_{h\to 0+}\frac{\norm{T_{t+h}u_{0}-T_{t}u_{0}}_{X}}{h}
    & = \limsup_{h\to
      0+}\frac{\norm{T_{\frac{t}{2}+h}(T_{\frac{t}{2}}u_{0})-T_{\frac{t}{2}}(T_{\frac{t}{2}}u_{0})}_{X}}{h}\\
    &\le 4\,L\, e^{\omega
      \frac{t}{2}}\,\frac{\norm{T_{t/2}u_{0}}_{X}}{t}\\
    &\le 2^{\delta+2} L\, M\,\,e^{\frac{1}{2}(\omega+\hat{\omega}) t}
      \frac{\norm{u_{0}}_{Y}^{\gamma}}{t^{\delta+1}}.
  \end{align*}
\end{proof}

Next, suppose that there is a partial ordering ``$\le$'' on $X$ such that
$(X,\le)$ is an ordered vector space. Then, we can state the following
theorem.

\begin{theorem}
  \label{thm:2}
  Let $(X,\le)$ be an ordered vector
  space, $C$ be a subset of $X$, and $\{T_{t}\}_{t\ge 0}$ be a family of mappings $T_{t} : C\to C$
  defined on a subset $C\subseteq X$ satisfying
  \begin{equation}\label{eq:23}
    \text{for every $u_{0}$, $\hat{u}_{0}\in C$ satisfying
      $u_{0}\le\hat{u}_{0}$, one has }T_{t}u_{0}\le
    T_{t}\hat{u}_{0}\text{ for all $t\ge 0$.}
  \end{equation}
 and
 \begin{equation}
  \label{eq:21}
  \lambda^{-1}\,T_{\lambda t}u _{0} =T_{t}[\lambda^{-1}u _{0}]\qquad\text{for all
    $\lambda>0$, $t\ge 0$
    and $u _{0}\in C$.}
\end{equation}
 Then for every $u_{0}\in C$ satisfying $u_{0}\ge
 0$, one has
 \begin{equation}
   \label{eq:24}
   \frac{T_{t+h}u_{0}-T_{t}u_{0}}{h}\le \frac{1}{t}T_{t}u_{0}\qquad\text{for every $t$, $h>0$.}
 \end{equation}
\end{theorem}

Before giving the proof of Theorem~\ref{thm:2}, we state the following
definition.

\begin{definition}
  If $(X,\le)$ is an order vector space then a family $\{T_{t}\}_{t\ge 0}$ of mappings $T_{t} : C\to C$
  defined on a subset $C\subseteq X$ is called \emph{order preserving}
  if $\{T_{t}\}_{t\ge 0}$ satisfies~\eqref{eq:23}.
\end{definition}

\begin{proof}[Proof of Theorem~\ref{thm:2}]
  Since $\left(1+\tfrac{h}{t}\right)^{-1}<1$, one has that
  $\left(1+\tfrac{h}{t}\right)^{-1}u_{0}\le u_{0}$. Then,
  by~\eqref{eq:16bis} for $f\equiv 0$ and~\eqref{eq:23}, one finds
  \begin{align*}
        T_{t+h}u_{0}-T_{t}u_{0}&=\left(1+\tfrac{h}{t}\right)
        T_{t}\left[\left(1+\tfrac{h}{t}\right)^{-1}u_{0}\right]-T_{t}u_{0}\\
   &=T_{t}\left[\left(1+\tfrac{h}{t}\right)^{-1}u_{0}\right]
                              -T_{t}u_{0}+\tfrac{h}{t}
     T_{t}\left[\left(1+\tfrac{h}{t}\right)^{-1}u_{0}\right]\\
    &\le \tfrac{h}{t} T_{t}\left[\left(1+\tfrac{h}{t}\right)^{-1}u_{0}\right]\\
    &\le \tfrac{h}{t} T_{t}u_{0},
  \end{align*}
from where one sees that~\eqref{eq:24} holds.
\end{proof}

By Theorem~\ref{thm:2}, if the derivative $\frac{\td }{\dt_{+}}T_{t}u_{0}$ exists
(in $X$) at $t>0$, then we can state the following.

\begin{corollary}
  \label{cor:2}
  Under the hypotheses of Theorem~\ref{thm:2}, suppose that for
  $u_{0}\in C$ satisfying $u_{0}\ge 0$, the right hand-side derivative
  $\frac{\td }{\dt_{+}}T_{t}u_{0}$ exists (in $X$) at $t>0$, then
 \begin{displaymath}
   \frac{\td T_{t}u_{0}}{\dt_{+}}\le \frac{1}{t}T_{t}u_{0}.
 \end{displaymath}
\end{corollary}

Further, we can conclude from Theorem~\ref{thm:2} the following
result.

\begin{corollary}
  \label{cor:3}
  In addition to the hypotheses of Theorem~\ref{thm:2}, suppose that
  there is a linear functional $\Lambda : X\to \R$ satisfying
  \begin{align}\label{eq:28}
    &\Lambda x\ge 0\qquad\text{ for every $x\in X$ satisfying $x\ge
      0$}\\ \label{eq:29}
    &\Lambda T_{t}u_{0}=\Lambda u_{0}\text{ for every $t\ge0$ and $u_{0}\in X$ satisfying $u_{0}\ge
      0$.}
  \end{align}
  Then, the following estimate holds for each $\nu\in \{+,-\}$,
  \begin{equation}
    \label{eq:31}
    \Lambda [T_{t+h}u_{0}-T_{t}u_{0}]^{\nu}\le \frac{h}{t}\,\Lambda x
    \qquad\text{for all $t$, $h>0$, $u_{0}\in C$ with $u_{0}\ge 0$.}
  \end{equation}
\end{corollary}

\begin{example}
  If $X=L^{q}(\Sigma,\mu)$ for some $\Sigma$-measure
  space $(\Sigma,\mu)$ and $1\le q\le \infty$, then an example
  for $\Lambda$ satisfying~\eqref{eq:28} and~\eqref{eq:29} is given by
  \begin{displaymath}
    \Lambda x=\int_{\Sigma}x \td\mu\qquad\text{for every $x\in X$.}
  \end{displaymath}
\end{example}

\begin{proof}[Proof of Corollary~\ref{cor:3}]
  Let $u_{0}\in C$ with $u_{0}\ge 0$, and $t$, $h>0$. Then we note
  first that  by~\eqref{eq:29},
  \begin{displaymath}
    0=\Lambda T_{t+h}u_{0}-\Lambda T_{t}u_{0}=\Lambda(T_{t+h}u_{0}-T_{t}u_{0})
  \end{displaymath}
  and since
  \begin{displaymath}
    \Lambda(T_{t+h}u_{0}-T_{t}u_{0})=\Lambda\left[T_{t+h}u_{0}-T_{t}u_{0}\right]^{+}
    -\Lambda\left[T_{t+h}u_{0}-T_{t}u_{0}\right]^{-},
  \end{displaymath}
  one has that
  \begin{equation}
    \label{eq:32}
    \Lambda\left[T_{t+h}u_{0}-T_{t}u_{0}\right]^{+}=
    \Lambda\left[T_{t+h}u_{0}-T_{t}u_{0}\right]^{-}.
  \end{equation}
  Further, by Theorem~\ref{thm:2} and since $T_{t}u_{0}\ge 0$, it
  follows from the definition of $[x]^{+}=\max\{x,0\}$, ($x\in \R$), that
  \begin{equation}
    \label{eq:30}
   \left[T_{t+h}u_{0}-T_{t}u_{0}\right]^{+}\le \frac{h}{t} T_{t}u_{0}.
  \end{equation}
  By the linearity of $\Lambda$ and by~\eqref{eq:28}, one has that $x\le y$ yields
  $\Lambda x\le \Lambda y$. Thus applying $\Lambda$ to~\eqref{eq:30}
  leads to~\eqref{eq:31} for $\nu=``+''$. Moreover, by~\eqref{eq:32},
  inequality~\eqref{eq:31} also holds for $\nu=``-''$. This completes
  the proof of this corollary.
\end{proof}

For the last result of this section, we consider the following differential inclusion
\begin{equation}
  \label{eq:10bisF}
  \frac{\td u}{\dt}+A(u(t))+ F(u(t))\ni 0\qquad\text{on $(0,+\infty)$,}
\end{equation}
for some operator $A\subseteq X\times X$ and a Lipschitz-continuous
mapping $F : X\to X$ with Lipschitz constant $\omega\ge 0$ and
satisfying $F(0)=0$. As for the differential
inclusion~\eqref{eq:10bis} and the case $f\equiv 0$, suppose, there is
a subset $C\subseteq X$ and a family $\{T_{t}\}_{t\ge 0}$ of mappings
$T_{t} : C\to C$ associated with $A$ through the relation that for
every given $u_{0}\in C$, the function $u$
defined by~\eqref{eq:22} is the unique solution of~\eqref{eq:10bisF}
with initial value $u(0)=u_{0}$. On the other hand, setting
\begin{equation}
  \label{eq:43}
  f(t):=-F(u(t)),\qquad\text{($t\ge0$),}
\end{equation}
one has that
\begin{equation}
  \label{eq:45}
  T_{t}(u_{0},f)=u(t)=T_{t}u_{0}\qquad\text{ for every $t\ge 0$, $u_{0}\in C$.}
\end{equation}
Thus, by Theorem~\ref{thm:1bis} for $T=+\infty$ we have the following estimates.

\begin{corollary}\label{cor:1bisF}
  Let $F : X\to X$ be a Lipschitz continuous mapping with
  Lipschitz-constant $\omega>0$ and satisfying $F(0)=0$. Suppose,
  there is a subset $C\subseteq X$, and a family $\{T_{t}\}_{t\ge 0}$
  of mappings $T_{t} : C\to C$ satisfying
  \begin{equation}
    \label{eq:50}
    \norm{T_{t}u_{0}}_{X}\le e^{\omega
      t}\,\norm{u_{0}}_{X}\qquad\text{for all $t\ge 0$, $u_{0}\in C$,}
  \end{equation}
  and in relation with~\eqref{eq:45}, suppose that $\{T_{t}\}_{t\ge
    0}$ satisfies~\eqref{eq:21bis} and~\eqref{eq:20bis} for $f$ given
  by~\eqref{eq:43}. Then for every $u_{0}\in C$, and $t$, $h>0$ such
  that $\abs{h}/t<1$, one has that
\begin{equation}
  \label{eq:47}
  \lnorm{\frac{T_{t+h}u_{0}-T_{t}u_{0}}{h}}_{X}
    \le \left[2 e^{L 2\omega \int_{0}^{t}e^{-\omega
        s}\ds}+\omega
    \int_{0}^{t}e^{L 2 \omega \int_{s}^{t}e^{-\omega
        r}\dr}\ds\right] \frac{e^{\omega t}L\norm{u_{0}}_{X}}{t}.
\end{equation}
Moreover, if the derivative
  $\frac{\td }{\dt}T_{t}u_{0}$ exists (in $X$) for a.e. $t>0$, then
  \begin{equation}
    \label{eq:48}
       \lnorm{\frac{\td T_{t}u_{0}}{\dt}}_{X}\le e^{\omega t}L\left[2e^{L \omega \int_{0}^{t} e^{-\omega s} s\,\ds}
  +\omega \int_{0}^{t} e^{L \omega \int_{s}^{t} e^{-\omega r} r\,\dr}\,\ds\right]\frac{\norm{u_{0}}_{X}}{t}
  \end{equation}
  for a.e. $t>0$.
\end{corollary}

For the proof of this corollary, we will employ the following version of
Gronwall's lemma.

\begin{lemma}
  \label{lem:1}
  Let $a\in L^{1}(0,T)$, $B : [0,T]\to \R$ be an absolutely continuous
  function, and $v\in L^{\infty}(0,T)$ satisfy
  \begin{displaymath}
    v(t)\le \int_{0}^{t}a(s)v(s)\,\ds+B(t)\qquad\text{for a.e. $t\in (0,T)$.}
  \end{displaymath}
  Then,
  \begin{displaymath}
    v(t)\le
    B(0)\,e^{\int_{0}^{t}a(s)\,\ds}+\int_{0}^{t}e^{\int_{s}^{t}a(r)\,\dr}\,B'(s)\,\ds
    \qquad\text{for a.e. $t\in (0,T)$.}
  \end{displaymath}
\end{lemma}

We now give the proof of Corollary~\ref{cor:1bisF}.

\begin{proof}
  Let $u_{0}\in C$, and $t$, $h>0$ such that $\abs{h}/t<1$. Then, by
  the hypotheses of this corollary, we are in the position to apply
  Theorem~\ref{thm:1bis} to $T_{t}(u_{0},f)$ for $f$ given
  by~\eqref{eq:43}. Then by~\eqref{eq:25bis}, one finds
  \begin{displaymath}
  \begin{split}
    \lnorm{\frac{T_{t+h}u_{0}-T_{t}u_{0}}{\tfrac{h}{t}}}_{X}
    &\le L\;e^{\omega\,t}\Big[2\norm{u_{0}}_{X}+\int_{0}^{t}e^{-\omega s}\norm{F(T_{s}u_{0})}_{X}\ds +\Big.\\
      &\hspace{1cm}\left.+ \left(1+\tfrac{h}{t}\right) \,
      \int_{0}^{t}e^{-\omega
        s}\lnorm{\frac{F(T_{s+\tfrac{h}{t}s}u_{0})-F(T_{s}u_{0})}{\frac{h}{t}}}_{X}\ds\right].
  \end{split}
\end{displaymath}
Since $F$ is globally Lipschitz continuous with constant
$\omega>0$, $F(0)=0$ and by~\eqref{eq:50}, it follows that
  \begin{displaymath}
  \begin{split}
    \lnorm{\frac{T_{t+h}u_{0}-T_{t}u_{0}}{\tfrac{h}{t}}}_{X}
    &\le L\;e^{\omega\,t}\Big[(2+\omega t)
      \norm{u_{0}}_{X}+\Big.\\
      &\hspace{1cm}\left.+ \left(1+\tfrac{h}{t}\right) \,\omega
      \int_{0}^{t}e^{-\omega
        s}\lnorm{\frac{T_{s+\tfrac{h}{t}s}u_{0}-T_{s}u_{0}}{\frac{h}{t}}}_{X}\ds\right].
  \end{split}
\end{displaymath}
Since $\abs{h}/t<1$,
\begin{equation}
\label{eq:46}
\begin{split}
  e^{-\omega\,t}\lnorm{\frac{T_{t+h}u_{0}-T_{t}u_{0}}{\tfrac{h}{t}}}_{X}
  & \le L\;\Big[(2+\omega t)
    \norm{u_{0}}_{X}+\Big.\\
    &\hspace{1cm}\left.+ 2 \,\omega \int_{0}^{t}e^{-\omega
      s}\lnorm{\frac{T_{s+\tfrac{h}{t}s}u_{0}-T_{s}u_{0}}{\frac{h}{t}}}_{X}\ds\right].
\end{split}
\end{equation}
Due to~\eqref{eq:46}, we can apply Gronwall's lemma to
  \begin{displaymath}
    B(t)=L (2+\omega t)\,\norm{u_{0}}_{X}\quad\text{ and }\quad
    a(t)=L 2\omega\, e^{-\omega t}.
  \end{displaymath}
  Then, one sees that~\eqref{eq:47} holds. Now, suppose that the
  derivative $\frac{\td }{\dt}T_{t}u_{0}$ exists (in $X$) for
  a.e. $t>0$, then by~\eqref{eq:40}, the Lipschitz continuity of $F$
  and by~\eqref{eq:50}, one one has that
  \begin{displaymath}
      e^{-\omega t}\, t\lnorm{\frac{\td T_{t}u_{0}}{\dt}}_{X} \le \left[2
  \norm{u_{0}}_{X}+ \omega \int_{0}^{t}e^{-\omega s} s
  \lnorm{\frac{\td T_{s}u_{0}}{\ds}}_{X}\ds+\omega t\,\norm{u_{0}}_{X}\right]
  \end{displaymath}
  for a.e. $t>0$. Now, applying Gronwall's lemma to
  \begin{displaymath}
    B(t)=L (2+\omega t)\,\norm{u_{0}}_{X}\quad\text{ and }\quad
    a(t)=L \omega\, e^{-\omega t} t,
  \end{displaymath}
  leads to~\eqref{eq:48}. This completes the proof of this corollary.
\end{proof}

%
%
%
%
%
%
%
%
%

\section{Accretive operators of homogeneous order zero}
\label{sec:semigroups}

Suppose $X$ is Banach space with norm $\norm{\cdot}_{X}$. Then, we
begin this section with the following definition.

\begin{definition}\label{def:quasi-accretive}
  For $\omega\in \R$, an operator $A$ on $X$ is called \emph{$\omega$-quasi
    $m$-accretive operator} on $X$ if $A$ is \emph{accretive}, that is, for every
  $(u,v)$, $(\hat{u},\hat{v})\in A$ and every $\lambda\ge 0$,
  \begin{displaymath}
    \norm{u-\hat{u}}_{X}\le \norm{u-\hat{u}+\lambda (\omega (u-\hat{u})+ v-\hat{v})}_{X}.
  \end{displaymath}
  and if for $A$ the \emph{range condition}
  \begin{equation}
    \label{eq:11}
    Rg(I+\lambda A)=X\qquad\text{for some (or equivalently, for all)
      $\lambda>0$, $\lambda\,\omega<1$,}
  \end{equation}
  holds.
\end{definition}

If $A$ is $\omega$-quasi $m$-accretive operator, then the classical
existence theorem~\cite[Theorem~6.5]{Benilanbook}
(cf~\cite[Corollary~4.2]{MR2582280}), for every
$u_{0}\in \overline{D(A)}^{\mbox{}_{X}}$ and $f\in L^{1}(0,T;X)$,
there is a unique \emph{mild} solution $u \in C([0,T];X)$
of~\eqref{eq:10bis}.

\begin{definition}\label{def:mild-solution}
  For given $u_{0}\in \overline{D(A)}^{\mbox{}_{X}}$ and $f\in
  L^{1}(0,T;X)$, a function $u\in C([0,T];X)$ is called a \emph{mild
    solution} of the inhomogeneous
differential inclusion~\eqref{eq:10bis} with initial value $u_{0}$
if $u(0)=u_{0}$ and for every $\varepsilon>0$,
there is a \emph{partition} $\tau_{\varepsilon} : 0=t_{0}<t_{1}<\cdots <
t_{N}=T$ and a \emph{step function}
\begin{displaymath}
  u_{\varepsilon,N}(t)=u_{0}\,\mathds{1}_{\{t=0\}}(t)+\sum_{i=1}^{N}u_{i}\,\mathds{1}_{(t_{i-1},t_{i}]}(t)
  \qquad\text{for every $t\in [0,T]$}
\end{displaymath}
satisfying
\begin{align*}
  &t_{i}-t_{i-1}<\varepsilon\qquad\text{ for all $i=1,\dots,N$,}\\
  &\sum_{N=1}^{N}\int_{t_{i-1}}^{t_{i}}\norm{f(t)-\overline{f}_{i}}\,\dt<\varepsilon\qquad
    \text{where
    $\overline{f}_{i}:=\frac{1}{t_{i}-t_{i-1}}\int_{t_{i-1}}^{t_{i}}f(t)\,\dt$,}\\
  & \frac{u_{i}-u_{i-1}}{t_{i}-t_{i-1}}+A u_{i}\ni \overline{f}_{i}\qquad\text{ for all $i=1,\dots,N$,}
\end{align*}
and
\begin{displaymath}
  \sup_{t\in [0,T]}\norm{u(t)-u_{\varepsilon,N}(t)}_{X}<\varepsilon.
\end{displaymath}
\end{definition}

In particular, if $A$ is $\omega$-quasi $m$-accretive, and if
for given $u_{0}\in \overline{D(A)}^{\mbox{}_{X}}$, $f\in
L^{1}(0,T;X)$, the function $u : [0,T]\to X$ is the unique mild solution of~\eqref{eq:10bis} with initial
value $u(0)=u_{0}$, then by~\eqref{eq:36} the family $\{T_{t}\}_{t=0}^{T}$ defines a \emph{semigroup}
of \emph{$\omega$-quasi contractions}
$T_{t} : \overline{D(A)}^{\mbox{}_{X}}\times L^{1}(0,T;X)\to \overline{D(A)}^{\mbox{}_{X}}$ for
$C=\overline{D(A)}^{\mbox{}_{X}}$; that is, $\{T_{t}\}_{t=0}^{T}$ satisfies
\begin{itemize}
  \item (\emph{semigroup property}) $T_{t+s}=T_{t}\circ T_{s}$ for
    every $t$, $s\in [0,T]$;
  \item (\emph{strong continuity}) for every $(u_{0},f)\in
    \overline{D(A)}^{\mbox{}_{X}}\times L^{1}(0,T;X)$, $t\mapsto
    T_{t}(u_{0},f)$ belongs to $C([0,T];X)$;
  \item (\emph{$\omega$-quasi
contractivity}) $T_{t}$ satisfies~\eqref{eq:20bis}
\end{itemize}
Furthermore, keeping $f\equiv 0$ and only varying $u_{0}\in
\overline{D(A)}^{\mbox{}_{X}}$, shows that by
\begin{equation}
  \tag{\ref{eq:55}}
  T_{t}u_{0}=T_{t}(u_{0},0)\qquad\text{for every $t\ge 0$.}
\end{equation}
defines a strongly continuous semigroup  $\{T_{t}\}_{t\ge 0}$ of
$\omega$-quasi contractions $T_{t} : \overline{D(A)}^{\mbox{}_{X}}\to \overline{D(A)}^{\mbox{}_{X}}$.
For the family $\{T_{t}\}_{t\ge 0}$ on
$\overline{D(A)}^{\mbox{}_{X}}$, the operator
\begin{displaymath}
  A_{0}:=\Bigg\{(u_{0},v)\in X\times X\Bigg\vert\;\lim_{h\downarrow
    0}\frac{T_{h}(u_{0},0)-u_{0}}{h}=v\text{ in $X$}\Bigg\}
\end{displaymath}
is an $\omega$-quasi accretive well-defined mapping
$A_{0} : D(A_{0})\to X$ and called the \emph{infinitessimal generator}
of $\{T_{t}\}_{t\ge 0}$. Under additional conditions on the geometry
of the Banach space $X$ (see Definition~\ref{def:Radon}), one has that
$A_{0}\subseteq A$.  Thus, we say (ignoring the abuse of details) that
both families $\{T_{t}\}_{t=0}^{T}$ on
$\overline{D(A)}^{\mbox{}_{X}}\times L^{1}(0,T;X)$ and
$\{T_{t}\}_{t\ge 0}$ on $\overline{D(A)}^{\mbox{}_{X}}$ are
\emph{generated by $-A$}.\medskip


In application, usually $X$ is given by the Lebesgue space
$L^{\infty}(\Sigma,\mu)$ (or $L^{r}(\Sigma,\mu)$ for $1\le r<\infty$)
and $Y$ is given by $L^{1}(\Sigma,\mu)$ (or $L^{r}(\Sigma,\mu)$ for
some $1\le q<r$) for some $\sigma$-finite measure space
$(\Sigma,\mu)$. Then, $L^{1}$-$L^{\infty}$-decay estimates are
intimately connected with abstract Sobolev inequalities satisfied by the
infinitesimal generator $-A$ of the semigroup $\{T_{t}\}_{t\ge
  0}$. For more details to the linear semigroup theory we refer to the
monograph~\cite{MR1218884} and to~\cite{CoulHau2017} for the nonlinear
semigroup theory.
%

Moreover (cf~\cite[Chapter~4.3]{Benilanbook}), for given
$u_{0}\in \overline{D(A)}^{\mbox{}_{X}}$ and any step function
$f=\sum_{i=1}^{N}f_{i}\,\mathds{1}_{(t_{i-1}, t_{i}]}\in
L^{1}(0,T;X)$, let $u : [0,T]\to X$ given by
\begin{equation}
  \label{eq:51}
  u(t)= u_{0}\,\mathds{1}_{\{t=0\}}(t)+\sum_{i=1}^{N}u_{i}(t) \mathds{1}_{(t_{i-1},t_{i}]}(t)
\end{equation}
is the unique mild solution
of~\eqref{eq:10bis}, where $u_{i}$ is the unique mild solution of
\begin{equation}
  \label{eq:52}
  \frac{\td u_{i}}{\dt}+A(u_{i}(t))\ni f_{i}\quad\text{ on
    $(t_{i-1},t_{i})$, and }\quad u_{i}(t_{i-1})=u_{i-1}(t_{i-1}).
\end{equation}
Then for every $i=1,\dots, N$, the semigroup $\{T_{t}\}_{t=0}^{T}$ is
obtained by the \emph{exponential formula}
\begin{equation}
  \label{eq:27}
  T_{t}(u(t_{i-1}),f_{i})=u_{i}(t)=\lim_{n\to\infty} \left[J_{\frac{t-t_{i-1}}{n}}^{A_{i}}\right]^{n}u(t_{i-1})\qquad\text{in
  $C([t_{i-1},t_{i}];X)$}
\end{equation}
for every $i=1,\dots,N$, where for $\mu>0$,
$J_{\mu}^{A_{i}}=(I+\mu A_{i})^{-1}$ is the \emph{resolvent operator}
of the operator $A_{i}$ and $A_i:= \{ (x, y - f_i) \ : \ (x,y) \in A \}$.\medskip

As for classical solutions, the fact that $A$ is homogeneous
  of order zero, is also reflected in the notion of mild solution and
  so in $\{T_{t}\}_{t\ge 0}$. This is shown in our next lemma.

\begin{lemma}\label{acrethomog} Let $A$ be a $\omega$-quasi
  $m$-accretive and $\{T_{t}\}_{t\ge 0}$ be the semigroup on
$\overline{D(A)}^{\mbox{}_{X}}\times L^{1}_{loc}([0,+\infty);X)$
generated by $-A$. If $A$ is homogeneous of
  order zero, then $\{T_{t}\}_{t\ge0}$ satisfies~\eqref{eq:21bis}
for every $(u_{0},f)\in \overline{D(A)}^{\mbox{}_{X}}\times L^{1}(0,T;X)$.
\end{lemma}

\begin{proof}
 For every $\mu>0$, $v\in X$, and $\lambda>0$, one has that
  \begin{displaymath}
    J_{\mu}^{A_{i}}\left[\lambda^{-1}v\right]=u\qquad\text{if and only if }\qquad
    u+\mu A_{i}u\ni \lambda^{-1}v,
  \end{displaymath}
  which if $A$ is \emph{homogeneous of
    order zero}, is equivalent to
  \begin{displaymath}
    \lambda u+\lambda\mu A_{i}(\lambda u)\ni v   \qquad\text{or }\qquad
    J_{\lambda\mu}^{A_{i}}v=\lambda u.
  \end{displaymath}
  Therefore,
  \begin{equation}
    \label{eq:26}
    \lambda^{-1}\, J^{A_{i}}_{\lambda\mu}v=
    J_{\mu}^{A_{i}}\left[\lambda^{-1}v\right]\qquad\text{for all $\lambda$,
      $\mu>0$, $v\in X$.}
  \end{equation}
 Now, for $u_{0}\in
 \overline{D(A)}^{\mbox{}_{X}}$ and a partition
 \begin{displaymath}
  \pi : 0=t_{0}<t_{1}<\cdots < t_{N}=T\quad\text{ of $[0,T]$}
\end{displaymath}
let $f=\sum_{i=1}^{N}f_{i}\mathds{1}_{(t_{i-1},t_{i}]}\in L^{1}(0,T;X)$
be a step function and $u$ be the unique mild solution
of~\eqref{eq:10bis} for $f$. Then $u$ is given by~\eqref{eq:51}, were
on each subinterval $(t_{i-1},t_{i}]$, $u_{i}$ is the unique mild
solution of~\eqref{eq:52}. For $t>0$, $n\in \N$, and
$\lambda\in (0,1]$, apply~\eqref{eq:26} to
\begin{displaymath}
  \mu=\frac{t}{n}\qquad\text{and}\qquad
  v=J^{A_{1}}_{\lambda\frac{t}{n}}[\lambda^{-1}u_{0}].
\end{displaymath}
Then,
\begin{displaymath}
  \left[J^{A_{1}}_{\frac{t}{n}}\right]^{2}[\lambda^{-1}u_{0}]=
  J_{\frac{t}{n}}^{A_{1}}\left[\lambda^{-1}J_{\lambda\frac{t}{n}}^{A_{1}}u_{0}\right]
  =\lambda^{-1}\left[J^{A_{1}}_{\lambda\frac{t}{n}}\right]^{2}u_{0}.
\end{displaymath}
Iterating this equation $n$-times, one finds that
\begin{equation}
  \label{eq:37}
  \lambda^{-1}\, \left[J^{A_{1}}_{\lambda \frac{t}{n}}\right]^{n}u_{0}=
\left[ J^{A_{1}}_{\frac{t}{n}}\right]^{n}\left[\lambda^{-1}u_{0}\right]
\end{equation}
and so, by~\eqref{eq:27} sending $n\to +\infty$ in the latter
equation, yields on the one site
\begin{displaymath}
  \lim_{n\to+\infty}
  \lambda^{-1}\, \left[J^{A_{1}}_{\lambda \frac{t}{n}}\right]^{n}u_{0}=\lambda^{-1}u_{1}(\lambda t) =\lambda^{-1}u(\lambda t)
\end{displaymath}
for every $t\in [0,\frac{t_{1}}{\lambda}]$, and on the other side
\begin{displaymath}
  \lim_{n\to+\infty}\left[
    J^{A_{1}}_{\frac{t}{n}}\right]^{n}\left[\lambda^{-1}u_{0}\right]=v(t)
\end{displaymath}
for every $t\in [0,\frac{t_{1}}{\lambda}]$, where $v$ is the unique
mild solution of~\eqref{eq:52} for $i=1$ on
$(0,\frac{t_{1}}{\lambda})$ with initial value
$v(0)=\lambda^{-1}u_{0}$. By uniqueness of the two limits, we have
thereby shown that
\begin{displaymath}
  \lambda^{-1}T_{\lambda t}(u_{0},f_{1})
 =
 T_{t}(\lambda^{-1}u_{0},f_{1}\mathds{1}_{(0,\frac{t_{1}}{\lambda}]})\qquad\text{for
   every $t\in
    \left[0,\frac{t_{1}}{\lambda}\right]$.}
\end{displaymath}
Similarly, for every $i=2,3,\dots, N$, replacing in~\eqref{eq:37} $u_{0}$ by
$u(t_{i-1})$ (where $u(t_{i-1})=u(\lambda \frac{t_{i-1}}{\lambda})=v(\frac{t_{i-1}}{\lambda})$),
$A_{1}$ by $A_{i}$, and $\frac{t}{n}$ by $\frac{t-\frac{t_{i-1}}{\lambda}}{n}$
gives
\begin{displaymath}
   \lambda^{-1}\, \left[J^{A_{i}}_{\lambda
       \frac{t-\frac{t_{i-1}}{\lambda}}{n}}\right]^{n}u(t_{i-1})=
\left[ J^{A_{i}}_{\frac{t-\frac{t_{i-1}}{\lambda}}{n}}\right]^{n}\left[\lambda^{-1}v(\frac{t_{i-1}}{\lambda})\right]
\end{displaymath}
and by sending  $n\to +\infty$, limit~\eqref{eq:27} leads one one side
to
\begin{displaymath}
  \lim_{n\to+\infty}\lambda^{-1}\, \left[J^{A_{i}}_{\lambda
   \frac{t-\frac{t_{i-1}}{\lambda}}{n}}\right]^{n}u(t_{i-1})=\lambda^{-1}u(\lambda t)
\end{displaymath}
and on the other side,
\begin{displaymath}
  \lim_{n\to+\infty}\left[
    J^{A_{i}}_{\frac{t-\frac{t_{i-1}}{\lambda}}{n}}\right]^{n}\left[\lambda^{-1}v(\frac{t_{i-1}}{\lambda})\right]=
  v(t)
\end{displaymath}
for every $t\in
\left[\frac{t_{i-1}}{\lambda},\frac{t_{i}}{\lambda}\right]$, where $v$
is the unique mild solution of~\eqref{eq:52} for $i$ on
$\left(\frac{t_{i-1}}{\lambda},\frac{t_{i}}{\lambda}\right)$ with initial value
$v(\frac{t_{i-1}}{\lambda})=\lambda^{-1}v(\frac{t_{i-1}}{\lambda})=\lambda^{-1}u(t_{i-1})$. Therefore,
and since $u$ is given by~\eqref{eq:51}, we have shown that
\begin{displaymath}
  \lambda^{-1}T_{\lambda
    t}(u(t_{i-1}),f_{i})=T_{t}(\lambda^{-1}u(t_{i-1}),f_{i}\mathds{1}_{(\frac{t_{i-1}}{\lambda},
    \frac{t_{i}}{\lambda}]})\qquad\text{for $t\in
    \left[\frac{t_{i-1}}{\lambda},\frac{t_{i}}{\lambda}\right]$.}
\end{displaymath}
Since for every step function $f$ on a partition $\pi$ of $[0,T]$, $u$
is given by~\eqref{eq:51}, we have thereby shown that~\eqref{eq:21bis}
holds if $f$ is a step function. Now, by~\eqref{eq:20bis}, an
approximation argument shows that if $A$ is homogeneous of order
zero, then the semigroup $\{T_{t}\}_{t\ge 0}$ on
$\overline{D(A)}^{\mbox{}_{X}}\times L^{1}(0,T;X)$ generated by $-A$
satisfies~\eqref{eq:21bis}.
 \end{proof}

 By the above Lemma and Theorem~\ref{thm:1bis}, we can now state the following.

\begin{corollary}\label{cor:1bis}
  For $\omega\in \R$, suppose $A$ is an $\omega$-quasi $m$-accretive
  operator on a Banach space $X$, and $A$ is homogeneous of order
  zero satisfying $0\in A0$. Then, for
  every $(u_{0},f)\in \overline{D(A)}^{\mbox{}_{X}}\times L^{1}(0,T;X)$, the semigroup
  $\{T_{t}\}_{t\ge 0}$ of mapping $T_{t} :
  \overline{D(A)}^{\mbox{}_{X}}\times L^{1}(0,T;X)\to
  \overline{D(A)}^{\mbox{}_{X}}$ generated by $-A$
  satisfies~\eqref{eq:25bis}-\eqref{eq:40}.
\end{corollary}

For having that regularity estimate~\eqref{eq:40}
(respectively,~\eqref{eq:18bis2}) is satisfied by the semigroup
$\{T_{t}\}_{t\ge 0}$, one requires that each mild solution $u$
of~\eqref{eq:10bis} (respectively, of~\eqref{eq:10}) is
\emph{differentiable} and a \emph{stronger} notion of solutions
of~\eqref{eq:10bis}. The next definition is taken
from~\cite[Definition~1.2]{Benilanbook}
(cf~\cite[Chapter~4]{MR2582280}).

\begin{definition}\label{def:strong-sols}
  A locally absolutely continuous function $u [0,T] : \to X$ is called a \emph{strong
      solution} of differential inclusion~\eqref{eq:10} if $u$ is
    differentiable a.e. on $(0,T)$, and for a.e. $t\in (0,T)$,
    $u(t)\in D(A)$ and $f(t)-\frac{\td u}{\dt}(t)\in A(u(t))$.
\end{definition}

The next characterization of strong solutions of~\eqref{eq:10}
highlights the important point of \emph{a.e. differentiability}.

\begin{proposition}[{\cite[Theorem~7.1]{Benilanbook}}]\label{propo:diff-mild}
  Let $X$ be a Banach space, $f\in L^{1}(0,T;X)$ and for $\omega\in \R$, $A$ be
  $\omega$-quasi $m$-accretive in $X$. Then $u$ is a strong solution
  of the differential inclusion~\eqref{eq:10bis} on $[0,T]$ if and only
  if $u$ is a mild solution on $[0,T]$ and $u$ is ``absolutely
  continuous'' on $[0,T]$ and differentiable a.e. on $(0,T)$.
\end{proposition}

Of course, every strong solution $u$ of~\eqref{eq:10bis} is a
  mild solution of~\eqref{eq:10bis}, absolutely continuous and
  differentiable a.e. on $[0,T]$. Moreover, the differential
  inclusion~\eqref{eq:10bis} admits mild and Lipschitz continuous
  solutions if $A$ is $\omega$-quasi $m$-accretive in $X$
  (cf~\cite[Lemma~7.8]{Benilanbook}). But absolutely continuous
  vector-valued functions $u : [0,T]\to X$ are not, in general,
  differentiable a.e. on $(0,T)$. However, if one assumes additional
  geometric properties on $X$, then the latter implication holds true. Our next
definition is taken from~\cite[Definition~7.6]{Benilanbook}
(cf~\cite[Chapter~1]{MR2798103}).\medskip

\begin{definition}\label{def:Radon}
  A Banach space $X$ is said to have the \emph{Radon-Nikod\'ym property}
  if every absolutely continuous function $F : [a,b]\to X$, ($a$,
  $b\in \R$, $a<b$), is differentiable almost everywhere on $(a,b)$.
\end{definition}

Known examples of Banach spaces $X$ admitting the Radon-Nikod\'ym property are:
\begin{itemize}
\item {\bfseries (Dunford-Pettis)} if $X=Y^{\ast}$ is separable, where
  $Y^{\ast}$ is the dual space of a Banach space $Y$;
\item if $X$ is \emph{reflexive}.
\end{itemize}

We emphasize that $X_{1}=L^{1}(\Sigma,\mu)$,
$X_{2}=L^{\infty}(\Sigma,\mu)$, or $X_{3}=C(\mathcal{M})$ for a
$\sigma$-finite measure space $(\Sigma,\mu)$, or respectively, for a
compact metric space $(\mathcal{M},d)$ don't have, in general, the
Radon-Nikod\'ym property (cf~\cite{MR2798103}). Thus, it is quite
surprising that there is a class of operators $A$ (namely, the class of
\emph{completely accretive operators}, see
Section~\ref{sec:completely-accretive} below), for which the
differential inclusion~\eqref{eq:10bis} nevertheless admits strong
solutions (with values in $L^{1}(\Sigma,\mu)$ or $L^{\infty}(\Sigma,\mu)$).\bigskip

Now, by Corollary~\ref{cor:1bis} and
Proposition~\ref{propo:diff-mild}, we can conclude the following
results. We emphasize that one crucial point in the statement of
Corollary~\ref{cor:RN} below is that due to the uniform
estimate~\eqref{eq:40}, one has that for all initial values
$u_{0}\in \overline{D(A)}^{\mbox{}_{X}}$, the unique strong solution
$u$ of~\eqref{eq:10bis} satisfying $u(0)=u_{0}$ is a strong solution,
and not only for $u_{0}\in D(A)$.

\begin{corollary}\label{cor:RN}
  For $\omega\in \R$, suppose $A$ is an $\omega$-quasi $m$-accretive
  operator on a Banach space $X$ admitting the Radon-Nikod\'ym property,
  and $\{T_{t}\}_{t\ge 0}$ is the semigroup on
  $\overline{D(A)}^{\mbox{}_{X}}\times L^{1}(0,T;X)$ generated by
  $-A$. If $A$ is homogeneous of order zero satisfying $0\in A0$,
  then for every $u_{0}\in \overline{D(A)}^{\mbox{}_{X}}$ and $f\in BV(0,T;X)$, the unique mild
  solution $u$ of~\eqref{eq:10bis} satisfying $u(0)=u_{0}$ is a strong solution and
  satisfies~\eqref{eq:40} for every $t>0$.
\end{corollary}

Now by Corollary~\ref{cor:1bisF} and
Proposition~\ref{propo:diff-mild}, we obtain the following result when
$A$ is perturbed by a Lipschitz mapping.

\begin{corollary}
  \label{cor:RN-Lipschitz-case}
 Suppose $X$ is a Banach space with the Radon-Nikod\'ym property,
  $F : X\to X$ be a Lipschitz continuous mapping with
  Lipschitz-constant $\omega>0$ satisfying $F(0)=0$, $A$ an
  $m$-accretive operator on $X$, and $\{T_{t}\}_{t\ge 0}$ is the
  semigroup on $\overline{D(A)}^{\mbox{}_{X}}$ generated by
  $-(A+F)$. If $A$ is homogeneous of order zero satisfying $0\in A0$,
  then~\eqref{eq:48} holds for every
  $u_{0}\in \overline{D(A)}^{\mbox{}_{X}}$ and a.e. $t>0$.
\end{corollary}

If the Banach space $X$ and its dual space $X^{\ast}$ are
\emph{uniformly convex}, then (cf~\cite[Theorem~4.6]{MR2582280}) for every
$u_{0}\in D(A)$, $f\in W^{1,1}(0,T;X)$, the mild solution
$u(t)=T_{t}(u_{0},f)$, ($t\ge 0$), of~\eqref{eq:10bis} is a strong
solution of~\eqref{eq:10bis}, $u$ is everywhere differentiable from
the right, $\frac{\td u}{\dt_{+}}$ is right continuous, and
\begin{displaymath}
  \frac{\td u}{\dt_{+}}(t)+(A-f(t))^{\! \circ}u(t)=0\qquad\text{for
    every $t\ge 0$,}
\end{displaymath}
where for every $t\in [0,T]$, $(A-f(t))^{\! \circ}$ denotes the \emph{minimal
  selection} of $A-f(t)$ defined by
\begin{displaymath}
  (A-f(t))^{\! \circ}:=\Big\{(u,v)\in
  A-f(t)\,\Big\vert\big. \norm{v}_{X}=\inf_{\hat{v}\in Au-f(t)}\norm{\hat{v}}_{X}\Big\}.
\end{displaymath}
Thus, under those assumptions on $X$ and by
Proposition~\ref{propo:diff-mild}, we can state the following three
corollaries. We begin by stating the inhomogeneous case.

\begin{corollary}
  Suppose $X$ and its dual space $X^{\ast}$ are
  uniformly convex, for $\omega\in \R$, $A$ is an
  $\omega$-quasi $m$-accretive operator on $X$, and
  $\{T_{t}\}_{t\ge 0}$ is the semigroup on
  $\overline{D(A)}^{\mbox{}_{X}}\times L^{1}(0,T;X)$ generated by $-A$. If $A$ is
  homogeneous of order zero satisfying $0\in A0$, then for every
  $u_{0}\in \overline{D(A)}^{\mbox{}_{X}}$ and $f\in W^{1,1}(0,T;X)$,
  \begin{displaymath}
    \begin{split}
      &\lnorm{(A-f(t))^{\! \circ}T_{t}(u_{0},f)}_{X}\\
      &\qquad\qquad \le \frac{e^{\omega\,t}}{t}\left[2 \norm{u_{0}}_{X}+
        \int_{0}^{t}e^{-\omega
          s}s\norm{f'(s)}_{X}\,\ds+\int_{0}^{t}e^{-\omega
          s}\norm{f(s)}_{X}\ds\right]
    \end{split}
  \end{displaymath}
  for every $t>0$.
\end{corollary}

The following corollary states the homogeneous case.

\begin{corollary}
  Suppose $X$ and its dual space $X^{\ast}$ are
  uniformly convex, for $\omega\in \R$, $A$ is an
  $\omega$-quasi $m$-accretive operator on $X$, and
  $\{T_{t}\}_{t\ge 0}$ is the semigroup on $\overline{D(A)}^{\mbox{}_{X}}$ generated by $-A$. If $A$ is
  homogeneous of order zero satisfying $0\in A0$, then
  \begin{displaymath}
  \lnorm{A^{\! \circ}T_{t}u_{0}}_{X}\le
  2 e^{\omega t}\frac{\norm{u_{0}}_{X}}{t}\qquad
  \text{ for every $t>0$ and $u_{0}\in \overline{D(A)}^{\mbox{}_{X}}$.}
  \end{displaymath}
\end{corollary}

The last corollary states the case when $A$ is perturbed by a
Lipschitz mapping. This follows from~\cite[Theorem~4.6]{MR2582280} and
Corollary~\ref{cor:1bisF}.

\begin{corollary}
  \label{cor:Lipschitz-continuous-case}
  Suppose $X$ and its dual space $X^{\ast}$ are uniformly convex,
  $F : X\to X$ be a Lipschitz continuous mapping with
  Lipschitz-constant $\omega>0$ satisfying $F(0)=0$, $A$ an
  $m$-accretive operator on $X$, and $\{T_{t}\}_{t\ge 0}$ is the
  semigroup on $\overline{D(A)}^{\mbox{}_{X}}$ generated by
  $-(A+F)$. If $A$ is homogeneous of order zero satisfying $0\in A0$,
  then for every $u_{0}\in \overline{D(A)}^{\mbox{}_{X}}$,
  \begin{displaymath}
       \lnorm{\frac{\td T_{t}u_{0}}{\dt}_{\! \!+}}_{X}\le e^{\omega t}\left[2e^{L \omega \int_{0}^{t} e^{-\omega s} s\,\ds}
  +\omega \int_{0}^{t} e^{L \omega \int_{s}^{t} e^{-\omega r} r\,\dr}\,\ds\right]\frac{\norm{u_{0}}_{X}}{t}
  \end{displaymath}
  for every $t>0$.
\end{corollary}

%
%

\section{Completely accretive operators of homogeneous
  order zero}
\label{sec:completely-accretive}

In \cite{MR1164641}, B\'enilan and Crandall introduced the celebrated
class of \emph{completely accretive} operators $A$ and showed that
there spaces without the Radon-Nikod\'ym property, but if $A$
  is homogeneous of order $\alpha>0$ with $\alpha\neq 1$, then the
  mild solutions of differential inclusion~\eqref{eq:10} involving
  $A$ are strong solutions. In this section we will see that this
also happen for completely accretive of homogeneous order
zero.\medskip

\subsection{General framework}

We begin by outlining our framework and then provide a brief
introduction to the class of completely accretive operators.\medskip

For the rest of this paper, suppose $(\Sigma, \mathcal{B}, \mu)$ is a $\sigma$-finite
measure space, and $M(\Sigma,\mu)$ the space of $\mu$-a.e. equivalent
classes of measurable functions $u : \Sigma\to \R$. For
$u\in M(\Sigma,\mu)$, we write $[u]^+$ to denote $\max\{u,0\}$ and
$[u]^-=-\min\{u,0\}$.  We denote by
$L^{q}(\Sigma,\mu)$, $1\le q\le \infty$, the corresponding standard
Lebesgue space with norm
\begin{displaymath}
  \norm{\cdot}_{q}=
  \begin{cases}
    \displaystyle\left(\int_{\Sigma}\abs{u}^{q}\,\textrm{d}\mu\right)^{1/q} &
    \text{if $1\le q<\infty$,}\\[7pt]
    \inf\Big\{k\in [0,+\infty]\;\Big\vert\;\abs{u}\le k\text{
      $\mu$-a.e. on $\Sigma$}\Big\}
    & \text{if $q=\infty$.}
  \end{cases}
\end{displaymath}
For $1\le q<\infty$, we
identify the dual space $(L^{q}(\Sigma,\mu))'$ with
$L^{q^{\mbox{}_{\prime}}}(\Sigma,\mu)$,
where $q^{\mbox{}_{\prime}}$ is the conjugate exponent of $q$ given by
$1=\tfrac{1}{q}+\tfrac{1}{q^{\mbox{}_{\prime}}}$.\medskip

Next, we first briefly recall the notion of \emph{Orlicz spaces}
(cf~\cite[Chapter 3]{MR1113700}). A continuous function
$\psi : [0,+\infty) \to [0,+\infty)$ is an \emph{$N$-function} if it
is convex, $\psi(s)=0$ if and only if $s=0$,
$\lim_{s\to 0+} \psi (s) / s = 0$, and
$\lim_{s\to\infty} \psi (s) / s = \infty$. Given an $N$-function
$\psi$, the \emph{Orlicz space}
  \begin{displaymath}
    L^\psi (\Sigma,\mu) := \Bigg\{ u \in M(\Sigma,\mu) \Bigg\vert\;
    \int_\Sigma
    \psi\left(\frac{|u|}{\alpha}\right)\; \td\mu <\infty \text{ for some }
    \alpha >0\Bigg\}
  \end{displaymath}
  and equipped with the \emph{Orlicz-Minkowski norm}
  \begin{displaymath}
    \norm{u}_{\psi} := \inf \Bigg\{ \alpha >0 \;\Bigg\vert\; \int_\Sigma \psi
    \left(\frac{|u|}{\alpha}\right)\; \td\mu\le 1 \Bigg\} .
  \end{displaymath}

With those preliminaries in mind, we are now in the position to recall the notation of
\emph{completely accretive} operators introduced in \cite{MR1164641} and
further developed to the \emph{$\omega$-quasi} case in~\cite{CoulHau2017}.\medskip

Let $J_{0}$ be the set given by
  \begin{displaymath}
    J_0 = \Big\{ j : \R \rightarrow
    [0,+\infty]\;\Big\vert\big. \text{$j$ is convex, lower
      semicontinuous, }j(0) = 0 \Big\}.
  \end{displaymath}
  Then, for every $u$, $v\in M(\Sigma,\mu)$, we write
  \begin{displaymath}
  u\ll v \quad \text{if and only if}\quad \int_{\Sigma} j(u)
  \,\td\mu \le \int_{\Sigma} j(v) \, \td\mu\quad\text{for all $j\in J_{0}$.}
\end{displaymath}

\begin{remark}
  Due to the interpolation result~\cite[Proposition~1.2]{MR1164641},
  for given $u$, $v\in M(\Sigma,\mu)$, the relation $u\ll v$ is
  equivalent to the two conditions
  \begin{displaymath}
    \begin{cases}
      \displaystyle\int_{\Sigma}[u-k]^{+}\,\td\mu&\le
      \int_{\Sigma}[v-k]^{+}\,\td\mu\qquad\text{for all $k>0$ and}\\[7pt]
       \displaystyle \int_{\Sigma}[u+k]^{-}\,\td\mu&\le
      \int_{\Sigma}[v+k]^{-}\,\td\mu\qquad\text{for all $k>0$.}
    \end{cases}
  \end{displaymath}
  Thus, the relation $\ll$ is closely related to the theory of
  rearrangement-invariant function spaces
  (cf~\cite{MR928802}). Another, useful characterization of relation
  $``\ll''$ is the following (cf~\cite[Remark~1.5]{MR1164641}): for
  $u$, $v\in M(\Sigma,\mu)$, $u\ll v$ if and only if $u^+\ll v^+$ and
  $u^-\ll v^-$.
\end{remark}

Further, the relation $\ll$ on $M(\Sigma,\mu)$ has the following
properties. We omit the easy proof of this proposition.

\begin{proposition}\label{propo:properties-of-ll}
  For every $u$, $v$, $w\in M(\Sigma,\mu)$, one has that
  \begin{enumerate}[(1)]
    \item\label{propo:properties-of-ll-claim0} $u^+ \ll u$, $u^{-}\ll -u$;
   \item\label{propo:properties-of-ll-claim1} $u\ll v$ if and only if $u^+\ll v^+$ and $u^-\ll v^-$;
   \item\label{propo:properties-of-ll-claim2} (\emph{positive homogeneity}) if $u\ll v$ then
    $\alpha u\ll \alpha v$ for all $\alpha>0$;
   \item\label{propo:properties-of-ll-claim3} (\emph{transitivity}) if $u\ll v$ and $v\ll w$ then $u\ll w$;
    \item\label{propo:properties-of-ll-claim4} if $u\ll v$ then
      $\abs{u}\ll \abs{v}$;
    \item\label{propo:properties-of-ll-claim5} (\emph{convexity}) for every $u\in M(\Sigma,\mu)$, the set
      $\{w\;\vert\,w\ll u \}$ is convex.
  \end{enumerate}
\end{proposition}

With these preliminaries in mind, we can now state the following
definitions.

\begin{definition}
  A mapping $S : D(S)\to M(\Sigma,\mu)$ with domain $D(S)\subseteq
  M(\Sigma,\mu)$ is called a \emph{complete contraction} if
  \begin{displaymath}
    Su-S\hat{u}\ll u-\hat{u}\qquad\text{for every $u$, $\hat{u}\in D(S)$.}
  \end{displaymath}
  More generally, for $L>0$, we call $S$ to be an
  \emph{L-complete contraction} if
  \begin{displaymath}
    L^{-1}Su-L^{-1}S\hat{u}\ll u-\hat{u}\qquad\text{for every $u$, $\hat{u}\in D(S)$,}
  \end{displaymath}
  or for $L=e^{\omega t}$ with $\omega\in \R$ and $t\ge 0$, $S$ is
  then also called an \emph{$\omega$-quasi complete contraction}.
\end{definition}

\begin{remark}
  \label{rem:1}
  Choosing $j_{q}(\cdot)=\abs{[\cdot]^{+}}^{q}\in \mathcal{J}_{0}$ if
  $1\le q<\infty$ and
  $j_{\infty}(\cdot)=[[\cdot]^{+}-k]^{+}\in \mathcal{J}_{0}$ for
  $k\ge 0$ large enough if $q=\infty$, or
  $j_{\psi,\alpha}(\cdot)=\psi(\frac{[\cdot]^{+}}{\alpha})$ for any
  $N$-function $\psi$ and $\alpha>0$ shows that for each $L$-complete
  contraction $S : D(S)\to M(\Sigma,\mu)$ with domain
  $D(S)\subseteq M(\Sigma,\mu)$, the mapping $L^{-1}S$ is order-preserving and
  contractive respectively with respect to the $L^{q}$-norm for all $1\le q\le \infty$, and
  the $L^{\psi}$-norm for any $N$-function $\psi$.
\end{remark}

Now, we can state the
definition of completely accretive operators.

\begin{definition}
  \label{def:completely-accretive-operators}
  An operator $A$ on $M(\Sigma,\mu)$ is called \emph{completely
    accretive} if for every $\lambda>0$, the resolvent operator
  $J_{\lambda}$ of $A$ is a complete contraction, or equivalently, if
  for every $(u_1,v_1)$, $(u_{2},v_{2}) \in A$ and $\lambda >0$, one
  has that
  \begin{displaymath}
    u_1 - u_2 \ll u_1 - u_2 + \lambda (v_1 - v_2).
  \end{displaymath}
  If $X$ is a linear subspace of $M(\Sigma,\mu)$ and $A$ an operator
  on $X$, then $A$ is \emph{$m$-completely accretive on $X$} if $A$ is
  completely accretive and satisfies the \emph{range
  condition}~\eqref{eq:11}.
Further, for $\omega\in \R$, an operator $A$ on a linear subspace
$X\subseteq M(\Sigma,\mu)$ is called \emph{$\omega$-quasi
  ($m$)-completely accretive} in $X$ if $A+\omega I$ is
($m$)-completely accretive in $X$.
\end{definition}

Before stating a useful characterization of completely accretive
operators, we first need to introducing the following function
spaces. Let
\begin{displaymath}
  L^{1+\infty}(\Sigma,\mu):= L^1(\Sigma,\mu) +
  L^{\infty}(\Sigma,\mu)\;\text{ and }\;
  L^{1\cap \infty}(\Sigma,\mu):=L^{1}(\Sigma,\mu)\cap L^{\infty}(\Sigma,\mu)
\end{displaymath}
be the \emph{sum} and the \emph{intersection space} of
$L^1(\Sigma,\mu)$ and $L^{\infty}(\Sigma,\mu)$, which respectively
equipped with the norms
\begin{align*}
  \norm{u}_{1+\infty}&:= \inf \left\{ \Vert u_1 \Vert_1 + \Vert
  u_2 \Vert_{\infty} \Big\vert u = u_1 + u_2, \ u_1 \in L^1(\Sigma,\mu), u_2
  \in L^{\infty}(\Sigma,\mu) \right\},\\
    \norm{u}_{1 \cap \infty}&:= \max\big\{ \norm{u}_{1}, \norm{u}_{\infty}\big\}
\end{align*}
are Banach spaces. In fact, $L^{1+\infty}(\Sigma,\mu)$ and and
$L^{1\cap \infty}(\Sigma,\mu)$ are respectively the largest and the
smallest of the rearrangement-invariant Banach function spaces
(cf~\cite[Chapter~3.1]{MR928802}). If $\mu(\Sigma)$ is finite, then
$L^{1+\infty}(\Sigma,\mu)=L^{1}(\Sigma,\mu)$ with equivalent norms, but
if $\mu(\Sigma)=\infty$ then $L^{1+\infty}(\Sigma,\mu)$ contains
$\bigcup_{1\le q\le\infty}L^{q}(\Sigma,\mu)$. Further, we will employ
the space
\begin{displaymath}
  L_0(\Sigma,\mu):= \left\{ u \in M(\Sigma,\mu) \;\Big\vert\,\Big.
    \int_{\Sigma} \big[\abs{u}
   - k\big]^+\,\td\mu < \infty \text{ for all $k > 0$} \right\},
\end{displaymath}
which equipped with the $L^{1+\infty}$-norm is a closed subspace of
$L^{1+\infty}(\Sigma,\mu)$. In fact, one has (cf~\cite{MR1164641})
that
\begin{math}
  L_0(\Sigma,\mu) = \overline{L^1(\Sigma,\mu) \cap
  L^{\infty}(\Sigma,\mu)}^{\mbox{}_{1+\infty}}.
\end{math}
Since for every $k\ge 0$,
$T_{k}(s):=[\abs{s}-k]^+$ is a Lipschitz mapping $T_{k} : \R\to\R$ and
by Chebyshev's inequality, one see that
$L^{q}(\Sigma,\mu)\hookrightarrow L_0(\Sigma,\mu)$ for every
$1\le q<\infty$ (and $q=\infty$ if $\mu(\Sigma)<+\infty$), and
$L^{\psi}(\Sigma,\mu)\hookrightarrow L_0(\Sigma,\mu)$ for every
$N$-function $\psi$.

\begin{proposition}[{\cite{CoulHau2017}}]
  \label{prop:completely-accretive}
   Let $P_{0}$ denote the set of all functions $T\in C^{\infty}(\R)$
   satisfying $0\le T'\le 1$,
    $T'$  is compactly supported, and $x=0$ is not contained in the
   support $\textrm{supp}(T)$ of $T$. Then for $\omega\in \R$,
   an operator $A \subseteq L_{0}(\Sigma,\mu)\times
  L_{0}(\Sigma,\mu)$ is $\omega$-quasi completely accretive if and only if
   \begin{displaymath}
     \int_{\Sigma}T(u-\hat{u})(v-\hat{v})\,\dmu+\omega \int_{\Sigma}T(u-\hat{u})(u-\hat{u})\,\dmu\ge 0
   \end{displaymath}
   for every $T\in P_{0}$ and every $(u,v)$, $(\hat{u},\hat{v})\in A$.
 \end{proposition}

\begin{remark}
  For convenience, we dote the unique extension of
  $\{T_{t}\}_{t\ge 0}$ on $L^\psi(\Sigma,\mu)$ or $L^1(\Sigma,\mu)$
  again by $\{T_{t}\}_{t\ge 0}$.
\end{remark}

\begin{definition}
  A Banach space $X\subseteq M(\Sigma,\mu)$ with norm
  $\norm{\cdot}_{X}$ is called \emph{normal} if the norm
  $\norm{\cdot}_{X}$ has the following property:
  \begin{displaymath}
    \begin{cases}
      & \text{for every $u\in X$, $v\in M(\Sigma,\mu)$ satisfying}\quad
        v\ll u,\\
      & \text{one has that}\quad v\in X\quad\text{ and }\quad\norm{v}_{X}\le
      \norm{u}_{X}.
    \end{cases}
  \end{displaymath}
\end{definition}

Typical examples of normal Banach spaces $X\subseteq M(\Sigma,\mu)$
are Orlicz-spaces $L^{\psi}(\Sigma,\mu)$ for every $N$-function
$\psi$, $L^{q}(\Sigma,\mu)$, ($1\le q\le \infty$),
$L^{1\cap \infty}(\Sigma,\mu)$, $L_0(\Sigma,\mu)$, and
$L^{1+\infty}(\Sigma,\mu)$.

\begin{remark}
  \label{rem:2}
  It is important to point out that if $X$ is a normal Banach space,
  then for every $u\in X$, one always have that $u^+$, $u^-$ and
  $\abs{u}\in X$. To see this, recall that
  by~\eqref{propo:properties-of-ll-claim0}
  Proposition~\ref{propo:properties-of-ll}, if $u\in X$, then $u^+\ll
  u$ and $u^-\ll -u$. Thus, $u^+$ and $u^-\in X$ and since
  $\abs{u}=u^+ +u^-$, one also has that $\abs{u}\in X$.
\end{remark}

The dual space
$(L_0(\Sigma,\mu))^{\mbox{}_{\prime}}$ of $L_0(\Sigma,\mu)$ is isometrically
isomorphic to $L^{1\cap \infty}(\Sigma,\mu)$. Thus, a sequence
$(u_{n})_{n\ge 1}$ in $L_0(\Sigma,\mu)$ is said to be \emph{weakly
  convergent} to $u$ in $L_0(\Sigma,\mu)$ if
\begin{displaymath}
  \langle v,u_{n}\rangle:=\int_{\Sigma} v\,u_{n}\,\td\mu\to
  \int_{\Sigma} v\,u\,\td\mu\qquad
  \text{for every $v\in L^{1\cap \infty}(\Sigma,\mu)$.}
\end{displaymath}
For the rest of this paper, we write $\sigma(L_{0},L^{1\cap\infty})$ to denote
the \emph{weak topology} on $L_0(\Sigma,\mu)$. For this weak topology,
we have the following compactness result.

\begin{proposition}[{\cite[Proposition~2.11]{MR1164641}}]
  \label{propo:compactness-in-L0}
  Let $u\in L_0(\Sigma,\mu)$. Then, the following statements hold.
  \begin{enumerate}[(1)]
   \item\label{propo:compactness-in-L0-claim1} The set
    \begin{math}
      \Big\{v\in M(\Sigma,\mu)\;\Big\vert\;v\ll u\Big\}
    \end{math}
    is $\sigma(L_{0},L^{1\cap\infty})$-sequently compact in
    $L_0(\Sigma,\mu)$;
  \item\label{propo:compactness-in-L0-claim2} Let
    $X\subseteq M(\Sigma,\mu)$ be a normal Banach space satisfying
    $X\subseteq L_0(\Sigma,\mu)$ and
     \begin{equation}
       \label{eq:7}
       \begin{cases}
         &\text{for every $u\in X$,
           $(u_{n})_{n\ge 1}\subseteq M(\Sigma,\mu)$ with $u_{n}\ll u$
           for all $n\ge 1$}\\[7pt]
         &\text{and $\displaystyle\lim_{n\to+\infty}u_{n}(x)=u(x)$ $\mu$-a.e. on
           $\Sigma$, yields}
         \quad \displaystyle\lim_{n\to +\infty}u_{n}=u\text{ in $X$.}
       \end{cases}
     \end{equation}
     Then for every $u\in X$ and sequence
     $(u_{n})_{n\ge 1}\subseteq M(\Sigma,\mu)$ satisfying
     \begin{displaymath}
       u_{n}\ll u\text{ for
     all $n\ge 1$}\quad\text{and}\quad \lim_{n\to+\infty}u_{n}=u
     \text{ $\sigma(L_{0},L^{1\cap\infty})$-weakly in $X$,}
   \end{displaymath}
   one has that
     \begin{displaymath}
       \lim_{n\to +\infty}u_{n}=u\qquad\text{ in $X$.}
     \end{displaymath}
  \end{enumerate}
\end{proposition}

Note, examples of normal Banach spaces $X\subseteq L_0(\Sigma,\mu)$
satisfying~\eqref{eq:7} are $X=L^{p}(\Sigma,\mu)$ for $1\le p<\infty$
and $L_{0}(\Sigma,\mu)$.\medskip

To complete this section we state the following Proposition
summarizing statements from~\cite[Proposition~2.9 \&
Proposition~2.10]{CoulHau2017}, which we will need in the sequel
(cf~\cite{MR1164641} for the case $\omega=0$).

\begin{proposition}\label{propo:properties-of-St}
  For $\omega\in \R$, let $A$ be $\omega$-quasi completely accretive in
  $L_{0}(\Sigma,\mu)$.
  \begin{enumerate}[(1.)]
  \item\label{propo:properties-of-St-claim1} If there is a $\lambda_{0}>0$ such that
    $\textrm{Rg}(I+\lambda A)$ is dense in $L_{0}(\Sigma,\mu)$, then
    for the closure $\overline{A}$ of $A$ in $L_{0}(\Sigma,\mu)$ and
    every normal Banach space with $X\subseteq L_0(\Sigma,\mu)$, the
    restriction $\overline{A}_{X}:=\overline{A}\cap(X\times X)$ of
    $A$ on $X$ is the unique $\omega$-quasi $m$-completely accretive extension of the
    part $A_{X}=A \cap(X\times X)$ of $A$ in $X$.

  \item For a given normal Banach space $X\subseteq L_0(\Sigma,\mu)$,
    and $\omega\in \R$, suppose $A$ is $\omega$-quasi $m$-completely accretive in $X$, and
    $\{T_{t}\}_{t\ge 0}$ be the semigroup generated by $-A$ on
    $\overline{D(A)}^{\mbox{}_{X}}$. Further, let $\{S_{t}\}_{t\ge 0}$
    be the semigroup generated by $-\overline{A}$, where
    $\overline{A}$ denotes the closure of $A$ in
    $\overline{X}^{\mbox{}_{L^{1+\infty}}}$. Then, the following
    statements hold.
    \begin{enumerate}[(a)]
    \item\label{propo:properties-of-St-claim1} The semigroup
      $\{S_{t}\}_{t\ge 0}$ is $\omega$-quasi completely contractive on
      $\overline{D(A)}^{\mbox{}_{L^{1+\infty}}}$, $T_{t}$ is the
      restriction of $S_{t}$ on $\overline{D(A)}^{\mbox{}_{X}}$,
      $S_{t}$ is the closure of $T_{t}$ in $L^{1+\infty}(\Sigma,\mu)$,
      and
      \begin{equation}
        \label{eq:exp-form-St}
        S_{t}u_{0}= L^{1+\infty}(\Sigma,\mu)-\lim_{n\to+\infty}
        \left(I+\frac{t}{n}A\right)^{-n}u_{0}\quad\text{for all $u_{0}\in \overline{D(A)}^{\mbox{}_{L^{1+\infty}}}\cap X$;}
      \end{equation}
    \item\label{propo:properties-of-St-claim5} If there exists
      $u\in L^{1\cap \infty}(\Sigma,\mu)$ such that the orbit
      $\{T_{t}u\,\vert\,t\ge 0\}$ is locally bounded on $\R_+$ with
      values in $L^{1\cap\,\infty}(\Sigma,\mu)$, then, for every
      $N$-function $\psi$, the semigroup $\{T_{t}\}_{t\ge 0}$ can be
      extrapolated to a strongly continuous, order-preserving
      semigroup of $\omega$-quasi contractions on
      $\overline{\overline{D(A)}^{\mbox{}_{X}}\cap
        L^{1\cap\,\infty}(\Sigma,\mu)}^{\mbox{}_{L^\psi}}$
      (respectively, on $\overline{\overline{D(A)}^{\mbox{}_{X}}\cap
        L^{1\cap\,\infty}(\Sigma,\mu)}^{\mbox{}_{L^1}}$), and to an
      order-preserving semigroup of $\omega$-quasi contractions on
      $\overline{\overline{D(A)}^{\mbox{}_{X}}\cap
        L^{1\cap\,\infty}(\Sigma,\mu)}^{\mbox{}_{L^\infty}}$. We
      denote each extension of $T_{t}$ on on those spaces again by
      $T_{t}$.

    \item\label{propo:properties-of-St-claim2} The restriction
      $A_{X}:=\overline{A}\cap (X\times X)$ of $\overline{A}$ on $X$
      is the unique $\omega$-quasi $m$-complete extension of $A$ in
      $X$; that is, $A=A_{X}$.
    \item\label{propo:properties-of-St-claim3} The operator $A$ is
      sequentially closed in $X\times X$ equipped with the
      relative\newline
      $(L_0(\Sigma,\mu)\times
      (X,\sigma(L_{0},L^{1\cap\infty})))$-topology.
    \item\label{propo:properties-of-St-claim3bis} The domain of $A$ is
      characterized by
      \begin{displaymath}
        \mbox{}\quad D(A)=\Big\{u\in \overline{D(A)}^{\mbox{}_{L^{1+\infty}}}\cap
        X\,\Big\vert\; \exists\;v\in X\text{ s.t.
        }e^{-\omega t}\frac{S_{t}u-u}{t}\ll v\text{ for small $t>0$}\Big\};
      \end{displaymath}
    \item\label{propo:properties-of-St-claim4} For every $u\in D(A)$,
      one has that
      \begin{equation}
        \label{eq:limit-infinit-gen-in-omega-A}
        \lim_{t\to 0+}\frac{S_{t}u-u}{t}=-A^{\!
          \circ}u\qquad\text{strongly in $L_{0}(\Sigma,\mu)$.}
      \end{equation}
    \end{enumerate}
  \end{enumerate}
\end{proposition}

%
%

\subsection{The subclass of homogeneous operators of order zero}
\label{subsec:homogeneous-zero}
As mentioned in Section~\ref{sec:semigroups}, the Banach spaces
$X_{1}=L^{1}(\Sigma,\mu)$ and $X_{2}=L^{\infty}(\Sigma,\mu)$ don't
have, the Radon-Niko\-d\'ym property. But for the class of
quasi $m$-completely accretive operators $A$ defined on a normal
Banach space $X\subseteq M(\Sigma,\mu)$, for semigroup
$\{T_{t}\}_{t\ge 0}$ generated by $-A$, the time-derivative
$\tfrac{\td T_{t}u_{0}}{\dt_{+}}$ exists in $X$ at every $t>0$ for
every $u_{0}\in \overline{D(A)}^{\mbox{}_{X}}$. This fact follows from
the following compactness result. Here, the partial ordering $``\!\!\!\le''$
is the standard one defined by $u\le v$ for $u$, $v\in M(\Sigma,\mu)$
if $u(x)\le v(x)$ for $\mu$-a.e. $x\in \Sigma$, and we use the symbol
$\hookrightarrow$ for indicating continuous embeddings.

\begin{lemma}\label{lem:compacteness-of-complete-semigroup}
  Let $X\subseteq L_{0}(\Sigma,\mu)$ be a normal Banach space
  satisfying~\eqref{eq:7}. For $\omega\in \R$, let
  $\{T_{t}\}_{t\ge 0}$ be a family of mappings $T_{t} : C\to C$
  defined on a subset $C\subseteq X$ of $\omega$-quasi complete
  contractions satisfying~\eqref{eq:21} and $T_{t}0=0$ for all
  $t\ge 0$. Then, for every $u_{0}\in C$ and $t>0$, the set
  \begin{equation}
    \label{eq:53}
    \Bigg\{\frac{T_{t+h}u_{0}-T_{t}u_{0}}{h}\,\Bigg\vert\Bigg.\,h\neq
    0, t+h>0\Bigg\}
  \end{equation}
  is $\sigma(L_{0},L^{1\cap\infty})$- weakly sequently compact in
  $L_{0}(\Sigma,\mu)$.
\end{lemma}

\begin{proof}
  Let $u_{0}\in C$, $t>0$, and $h\neq 0$ such that $t+h>0$. Then by taking
  $\lambda=1+\frac{h}{t}$ in~\eqref{eq:21}, one sees that
  \begin{align*}
    \abs{T_{t+h}u_{0}-T_{t}u_{0}}&= \abs{\lambda
                             \,T_{t}\left[\lambda^{-1}u_{0}\right]-T_{t}u_{0}}\\
    & \le     \lambda\,\labs{T_{t}\left[\lambda^{-1}u_{0}\right]-T_{t}u_{0}}+\abs{1-\lambda}\,
    \abs{T_{t}u_{0}}.
  \end{align*}
  Since $T_{t}$ is an $\omega$-quasi complete contraction,
  by~\eqref{propo:properties-of-ll-claim2} of
  Proposition~\ref{propo:properties-of-ll}, and since $T_{t}0=0$,
  ($t\ge 0$), one has that
  \begin{displaymath}
    \lambda\,e^{-\omega t}\labs{T_{t}\left[\lambda^{-1}u_{0}\right]-T_{t}u_{0}}\ll
    \abs{1-\lambda}\,\abs{u_{0}}
  \end{displaymath}
  and
  \begin{displaymath}
    \abs{1-\lambda}\, e^{-\omega t}
    \abs{T_{t}u_{0}}\ll \abs{1-\lambda}\,
    \abs{u_{0}}.
  \end{displaymath}
  Since the set $\{w\,\vert\,w\ll \abs{1-\lambda}\,
  \abs{u_{0}}\}$ is convex (see~\eqref{propo:properties-of-ll-claim5} of
  Proposition~\ref{propo:properties-of-ll}), we can conclude that
  \begin{displaymath}
   \frac{1}{2} e^{-\omega t}\abs{T_{t+h}u_{0}-T_{t}u_{0}}\ll\, \abs{1-\lambda}\,
    \abs{u_{0}}=\frac{\abs{h}}{t}\,\abs{u_{0}}
  \end{displaymath}
  and hence, by~\eqref{propo:properties-of-ll-claim2} of
  Proposition~\ref{propo:properties-of-ll},
  \begin{equation}
    \label{eq:5}
    \frac{\abs{T_{t+h}u_{0}-T_{t}u_{0}}}{\abs{h}}\ll\,2\, e^{\omega t}\,\frac{\abs{u_{0}}}{t}.
  \end{equation}
  Since for every $u\in M(\Sigma,\mu)$, one always has that
  $u^+\ll\abs{u}$, the transitivity of $``\ll''$ (see~\eqref{propo:properties-of-ll-claim3} of
  Proposition~\ref{propo:properties-of-ll}) implies for
  \begin{displaymath}
   f_{h}:=\frac{T_{t+h}u_{0}-T_{t}u_{0}}{\abs{h}},\qquad\text{one has
     that}\qquad f_{h}^{+}\ll\,2\, e^{\omega t}\,\frac{\abs{u_{0}}}{t}
  \end{displaymath}
  Therefore, by~\eqref{propo:compactness-in-L0-claim1} of
  Proposition~\ref{propo:compactness-in-L0}, the two sets
  $\{f_{h}^{+}\vert\,h\neq 0, t+h>0\}$ and
  $\{\abs{f_{h}}\vert\,h\neq 0, t+h>0\}$ are
  $\sigma(L_{0},L^{1\cap\infty})$- weakly sequently compact in
  $L_{0}(\Sigma,\mu)$, and since $f_{h}^{-}=\abs{f_{h}}-f_{h}^{+}$ and
  $f_{h}=f_{h}^{+}-f_{h}^{-}$, we have thereby shown that the claim of
  this lemma holds.
\end{proof}

With these preliminaries in mind, we can now state the regularization
effect of the semigroup $\{T_{t}\}_{t\ge 0}$ generated by a
$\omega$-quasi $m$-completely accretive operator of homogeneous order
zero.

\begin{theorem}\label{thm:regularity-of-complete-acc-homogen-operators}
  Let $X\subseteq L_0(\Sigma,\mu)$ be a normal Banach space
  satisfying~\eqref{eq:7}. For $\omega\in \R$, let $A$ be $\omega$-quasi $m$-completely accretive
  in $X$, and $\{T_{t}\}_{t\ge 0}$ be the semigroup generated by $-A$
  on $\overline{D(A)}^{\mbox{}_{X}}$. If $(0,0)\in A$ and
  $A$ is homogeneous of order zero, then for every
  $u_{0}\in \overline{D(A)}^{\mbox{}_{X}}$ and $t>0$,
  $\frac{\td T_{t}u_{0}}{\dt}$ exists in $X$ and
  \begin{equation}
    \label{eq:3}
    \abs{A^{\! \circ}T_{t}u_{0}}\le 2 e^{\omega
      t}\,\frac{\abs{u_{0}}}{t}\qquad\text{$\mu$-a.e. on $\Sigma$.}
  \end{equation}
  In particular,
  \begin{equation}
    \label{eq:6}
    \lnorm{\frac{\td T_{t}u_{0}}{\dt}}\le 2 e^{\omega
      t}\,\frac{\norm{u_{0}}}{t}
   \qquad\text{for every $t>0$,}
  \end{equation}
  and
  \begin{equation}
    \label{eq:4}
    \frac{\td T_{t}u_{0}}{\dt}\le
    \frac{T_{t}u_{0}}{t}\qquad\text{$\mu$-a.e. on $\Sigma$ for every
      $t>0$ if $u_{0}\ge 0$,}
  \end{equation}
  for every $u_{0}\in \overline{D(A)}^{\mbox{}_{X}}$ (then
  $\norm{\cdot}$ denotes the norm on $X$), respectively, for every
  $u_{0}\in \overline{\overline{D(A)}^{\mbox{}_{X}}\cap
    L^{1\cap\infty}(\Sigma,\mu)}^{\mbox{}_{L^{\psi}}}$ (then
  $\norm{\cdot}$ is the $L^{\psi}$-norm) for every $N$-function $\psi$ or for
  every $1\le \psi\equiv p<\infty$, and
  for every $u_{0}\in \overline{D(A)}^{\mbox{}_{X}}\cap L^{\infty}(\Sigma,\mu)$
  (where then $\norm{\cdot}$ is the $L^{\infty}$-norm).
\end{theorem}

\begin{proof}
  Let $u_{0}\in \overline{D(A)}^{\mbox{}_{X}}$, $t>0$, and
  $(h_{n})_{n\ge 1}\subseteq \R$ be a zero sequence such that
  $t+h_{n}>0$ for all $n\ge 1$. Due to Lemma
  \ref{acrethomog}, we can apply
  Lemma~\ref{lem:compacteness-of-complete-semigroup}. Thus, there is a
  $z\in L_{0}(\Sigma,\mu)$ and a subsequence
  $(h_{k_{n}})_{n\ge 1}$ of $(h_{n})_{n\ge 1}$ such that
  \begin{equation}
    \label{eq:1}
    \lim_{n\to
      +\infty}\frac{T_{t+h_{k_{n}}}u_{0}-T_{t}u_{0}}{h_{k_{n}}}=z\qquad\text{weakly
      in $L_{0}(\Sigma,\mu)$.}
  \end{equation}
  Moreover, by~\eqref{propo:properties-of-St-claim3bis} of
  Proposition~\ref{propo:properties-of-St}, 
  one has that $(T_{t}u_{0},-z)\in A$. Thus
  \eqref{propo:properties-of-St-claim4} of the same
  proposition~\ref{propo:properties-of-St} yields that
  $z=-A^{\!\circ}T_{t}u_{0}$ and 
  \begin{equation}
    \label{eq:2-new}
    \lim_{n\to 0}\frac{T_{t+h_{k_{n}}}u_{0}-T_{t}u_{0}}{h_{k_{n}}}=-A^{\!\circ}u_{0}
    \qquad\text{strongly in $L_{0}(\Sigma,\mu)$.}
  \end{equation}
  After possibly passing to another subsequence, we have that
  limit~\eqref{eq:2-new} holds also $\mu$-a.e. on $\Sigma$. Since
  $2 e^{-\omega t}\,\frac{\abs{u_{0}}}{t}\in X$ and
  $X\subseteq L_{0}(\Sigma,\mu)$,
  \eqref{propo:compactness-in-L0-claim2} of
  Proposition~\ref{propo:compactness-in-L0} implies that
  \begin{equation}
    \label{eq:8}
    \lim_{h\to 0}\frac{T_{t+h}u_{0}-T_{t}u_{0}}{h}=-A^{\!
      \circ}T_{t}u_{0}\qquad\text{exists in
    $X$ and $\mu$-a.e. on $\Sigma$.}
  \end{equation}
  Thus and since by~\eqref{eq:5},
  \begin{displaymath}
    \frac{\abs{T_{t+h_{k_{n}}}u_{0}-T_{t}u_{0}}}{\abs{h_{k_{n}}}}\le\,2\,
    e^{-\omega t}\,\frac{\abs{u_{0}}}{t}
    \qquad\text{for all $n\ge 1$,}
  \end{displaymath}
  sending $n\to +\infty$ in the last inequality,
  gives~\eqref{eq:3}. In particular, by Corollary~\ref{thm:1}, one has
  that~\eqref{eq:6} holds for the norm $\norm{\cdot }_{X}$ on $X$ and
  by Theorem~\ref{thm:2} that~\eqref{eq:4} holds. Moreover, we have
  that $-A^{\! \circ}T_{t}u_{0}=\frac{\td T_{t}u_{0}}{\dt_{\!+}}$
  $\mu$-a.e. on $\Sigma$ for every $t>0$. Thus, sending $h\to0+$
  in~\eqref{eq:5} shows that~\eqref{eq:3} holds. Further, by
  the $\mu$-a.e.-limit~\eqref{eq:8}, applying Fatou's
  lemma to~\eqref{eq:5} yields that~\eqref{eq:6} holds for the
  $L^{\psi}$-norm for every $N$-function $\psi$ and the $L^{p}$-norm
  $1\le p<\infty$. Since~\eqref{eq:6} holds for all $p<\infty$,
  sending $1\le p\to+\infty$ completes the proof of this theorem.
\end{proof}

\begin{remark}[{\bfseries Open problem}]
  We emphasize that the crucial point in the previous proof is that
  due to the zero-order homogeneity of $A$, the set~\eqref{eq:53} is
  $\sigma(L_{0},L^{1\cap\infty})$- weakly sequently compact in
  $L_{0}(\Sigma,\mu)$ and hence, for every
  $t>0$, $T_{t}u_{0}\in D(A)$ and
  \begin{displaymath}
    \frac{\td T_{t}u_{0}}{\dt_{\!+}}=\lim_{h\to0+}\frac{T_{t+h}u_{0}-T_{t}u_{0}}{h}=-A^{\!
      \circ}T_{t}u_{0}\qquad\text{exists in $X$.}
\end{displaymath}

  We believe that this remains true if the infinitesimal generator of the semigroup $\{T_{t}\}_{t\ge0}$ is
  of the form $A+F$ where $A$ is homogeneous of order zero and $F$ is
  Lipschitz-continuous. But so far, we are not able to show this result.
\end{remark}

As a final result of this section, we state the following decay
estimates for semigroups generated by the perturbed operator
$A+F$. Here, we write $L^{q_{0}\cap\infty}$ for the intersection space
$L^{q_{0}}(\Sigma,\mu)\to L^{\infty}(\Sigma,\mu)$.

\begin{theorem}
  \label{thm:Lipschitz}
  Let $F : M(\Sigma,\mu)\to M(\Sigma,\mu)$ be a mapping such that for
  every $N$-function $\psi$ and for $\psi\equiv 1$ and
  $\psi\equiv +\infty$, the restriction
  $F_{\vert L^{\psi}}: L^{\psi}(\Sigma,\mu)\to L^{\psi}(\Sigma,\mu)$
  is Lipschitz continuous with constant Lipschitz $\omega>0$ and
  $F(0)=0$.
  Let $A$ be an $m$-completely accretive operator on normal Banach
  space $X\subseteq L_{0}(\Sigma,\mu)$, and
  $\{T_{t}\}_{t\ge 0}$ the semigroup generated by $-(A+F)$ on
  $\overline{D(A)}^{\mbox{}_{X}}$. If $(0,0)\in A$ and $A$ is
  homogeneous of order zero, then
  \begin{equation}
    \label{eq:6-F-Lq}
    \lnorm{\frac{\td T_{t}u_{0}}{\dt _{\!+}}}_{\psi}\le \left[2 e^{\omega \int_{0}^{t}e^{-\omega
          s}s\ds}+\omega\int_{0}^{t}e^{\omega\int_{s}^{t}e^{-\omega r}r\dr}\ds\right]  \frac{e^{\omega
      t}\,\norm{u_{0}}_{\psi}}{t}
  \end{equation}
  for every $t>0$, and 
  every $u_{0}\in \overline{\overline{D(A)}^{\mbox{}_{X}}\cap
    L^{1\cap\infty}(\Sigma,\mu)}^{\mbox{}_{L^{\psi}}}$ for every
  $N$-function $\psi$, and every $1\le \psi\equiv p<\infty$, and
  for every $u_{0}\in \overline{D(A)}^{\mbox{}_{X}}\cap
  L^{q_{0}\cap \infty}(\Sigma,\mu)$ for $q=\infty$
  (where then $\norm{\cdot}_{\psi}$ is the $L^{\infty}$-norm).
  %
\end{theorem}

\begin{proof}
  Since $(0,0)\in (A+F)$, $T_{t}0=0$ for all $t\ge 0$. Thus,
  $u\equiv0\in L^{1\cap\infty}(\Sigma,\mu)$ such that
  $\{T_{t}u\vert\,t\ge0\}$ is locally bounded in $\R_{+}$. Thus, by
  Proposition~\ref{propo:properties-of-St}, for every $N$-function
  $\psi$ (respectively, for $\psi\equiv 1$ and $\psi\equiv\infty$)
  each $T_{t}$ admits a unique extension (which we denote again by
  $T_{t}$) of an $\omega$-quasi contractions on
  $\overline{\overline{D(A)}^{\mbox{}_{X}}\cap
    L^{1\cap\infty}(\Sigma,\mu)}^{\mbox{}_{L^{\psi}}}$ with respect to
  the $L^{\psi}$-norm. In addition, the family $\{T_{t}\}_{t\ge 0}$
  remains a semigroup satisfying~\eqref{eq:50} and in relation
  with~\eqref{eq:45}, $\{T_{t}\}_{t\ge 0}$ satisfies~\eqref{eq:21bis}
  and~\eqref{eq:20bis} for $f$ given by~\eqref{eq:43}. Further, for
  $1<\psi\equiv q<\infty$, $L^{q}(\Sigma,\mu)$ and its dual space
  $L^{q^{\mbox{}_{\prime}}}(\Sigma,\mu)$ are uniformly
  convex. Therefore, by Corollary~\ref{cor:1bisF} and
  Proposition~\ref{propo:diff-mild}, for every
  $u_{0}\in \overline{\overline{D(A)}^{\mbox{}_{X}}\cap
    L^{1\cap\infty}(\Sigma,\mu)}^{\mbox{}_{L^{q}}}$, for every $t>0$,
  $\frac{\td T_{t}u_{0}}{\dt}_{\!\! +}$ exists in $L^{q}(\Sigma,\mu)$,
  \begin{equation}
    \label{eq:54}
    \frac{\td T_{t}u_{0}}{\dt_{\!+}}=\lim_{h\to0+}\frac{T_{t+h}u_{0}-T_{t}u_{0}}{h}\quad\text{exists
    $\mu$-a.e. on $\Sigma$,}
  \end{equation}
  and~\eqref{eq:6-F-Lq} holds. Moreover, by Corollary~\ref{cor:1bisF},
  one has that \eqref{eq:47} holds for the $L^{\psi}$-norm and every
  $u_{0}\in \overline{\overline{D(A)}^{\mbox{}_X}\cap
    L^{1\cap\infty}(\Sigma,\mu)}^{\mbox{}_{L^{\psi}}}$ and every
  $N$-function $\psi$, respectively for the $L^{1}$-norm and every
  $u_{0}\in \overline{\overline{D(A)}^{\mbox{}_{X}}\cap
    L^{1\cap\infty}(\Sigma,\mu)}^{\mbox{}_{L^{1}}}$. Thus and
  by~\eqref{eq:54}, sending $h\to0+$ in~\eqref{eq:47} one obtains
  that~\eqref{eq:6-F-Lq} holds for all $N$-function $\psi$ and $q=1$.

  Next, let
  $u_{0}\in\overline{D(A)}^{\mbox{}_{X}}\cap
  L^{q_{0}\cap\infty}(\Sigma,\mu)$ for some $1\le q_{0}<+\infty$ and
  $t>0$. We assume $\norm{\frac{\td T_{t}u_{0}}{\dt}_{\!\! +}}_{\infty}>0$
  (otherwise, there is nothing to show). Then, for every
  $s\in (0,\norm{\frac{\td T_{t}u_{0}}{\dt}_{\!\! +}}_{\infty})$ and every
  $q_{0}\le q<\infty$, Chebyshev's inequality yields
  \begin{displaymath}
    \mu\left(\Bigg\{\labs{\frac{\td T_{t}u_{0}}{\dt}_{\!\! +}}\ge s
        \Bigg\}\right)^{1/q}\le \frac{\lnorm{\frac{\td T_{t}u_{0}}{\dt}_{\!\! +}}_{q}}{s}
  \end{displaymath}
  and so, by~\eqref{eq:6-F-Lq},
  \begin{displaymath}
    s\,\mu\left(\Bigg\{\labs{\frac{\td T_{t}u_{0}}{\dt}_{\!\! +}}\ge s
        \Bigg\}\right)^{1/q}\le \left[2 e^{\omega \int_{0}^{t}e^{-\omega
          s}s\ds}+\omega\int_{0}^{t}e^{\omega\int_{s}^{t}e^{-\omega
          r}r\dr}\ds\right]
    \frac{e^{\omega t}\,\norm{u_{0}}_{q}}{t}.
  \end{displaymath}
 Thus and since $\lim_{q\to
    \infty}\norm{u_{0}}_{q}=\norm{u_{0}}_{\infty}$, sending
  $q\to+\infty$ in the last inequality, yields
\begin{displaymath}
    s \le \left[2 e^{\omega \int_{0}^{t}e^{-\omega
          s}s\ds}+\omega\int_{0}^{t}e^{\omega\int_{s}^{t}e^{-\omega r}r\dr}\ds\right]  \frac{e^{\omega
      t}\,\norm{u_{0}}_{\infty}}{t}
  \end{displaymath}
  and since
  $s\in (0,\lnorm{\frac{\td T_{t}u_{0}}{\dt}_{\!\!  +}}_{\infty})$ was
  arbitrary, we have thereby shown that~\eqref{eq:6-F-Lq} also holds
  for $q=\infty$.
\end{proof}

%
%

\section{Application}
\label{sec:application}

Throughout this section, let $\Sigma$ be an open set of $\R^{d}$ and
the Lebesgue space $L^{q}(\Sigma)$ is equipped with the classical
Lebesgue measure. Suppose $f : \Sigma\times\R\to \R$ is a
Lipschitz-continuous \emph{Carath\'eodory} function, that is, $f$
satisfies the following three properties:
\begin{align}
  \label{property:one-Charatheodory}
\bullet\quad & f(\cdot, u) : \Sigma\to \R \text{ is measurable on $\Sigma$ for every
$u\in \R$,}\\
  \label{property:two-Charatheodory}
\bullet\quad & f(x,0)=0\text{ for a.e. $x\in \Sigma$, and}\\
 \notag
\bullet\quad & \text{there is a constant $\omega\ge 0$ such that }\\ \label{eq:2}
    \begin{split}
    \qquad\abs{f(x,u)-f(x,\hat{u})}\le \omega\,\abs{u-\hat{u}}\quad\,
  \text{for all $u$, $\hat{u}\in \R$, a.e. $x\in \Sigma$.}
    \end{split}
\end{align}
Then, for every $1\le q\le \infty$, $F : L^{q}(\Sigma)\to L^{q}(\Sigma)$ defined by
\begin{displaymath}
  F(u)(x):=f(x,u(x))\qquad\text{ for every $u\in L^{q}(\Sigma)$}
\end{displaymath}
is the associated \emph{Nemytskii operator} on $L^{q}(\Sigma)$. Moreover,
by~\eqref{eq:2}, $F$ is globally Lipschitz continuous on
$L^{q}(\Sigma)$ with constant $\omega>0$ and $F(0)(x)=0$ for
a.e. $x\in \Sigma$.

\subsection{Decay estimates of the total variational flow}

In this subsection, we consider the perturbed total variational flow operator
($1$-Laplace operator) given by
\begin{displaymath}
  Au:=-\Delta_{1}u +f(x,u)\qquad\text{with}\qquad
  \Delta_{1}u=\textrm{div}\left(\frac{D u}{\abs{D u}}\right),
\end{displaymath}
equipped with either Neumann boundary conditions on a
bounded domain $\Sigma$ in $\R^{d}$, $d\ge 1$.\medskip

Here, we use the following notation. A function
$u \in L^1(\Sigma)$ is said to be a \emph{function of bounded
  variation in $\Sigma$}, if the distributional partial derivatives
$D_{1}u:=\tfrac{\partial u}{\partial x_{1}}$, $\dots,$
$D_{d}u:=\tfrac{\partial u}{\partial x_{d}}$ are finite Radon measures in $\Sigma$, that is, if
\begin{displaymath}
  \int_{\Omega} u\,D_{i}\varphi\,\dx= -
  \int_{\Omega}\varphi\,\textrm{d} D_{i}u
\end{displaymath}
for all $\varphi\in C^{\infty}_{c}(\Sigma)$, $i=1, \dots, d$. The
linear vector space of functions $u \in L^1(\Sigma)$ of bounded
variation in $\Sigma$ is denoted by $BV(\Sigma)$. Further, we set
$Du=(D_{1}u,\dots,D_{d}u)$ for the \emph{distributional gradient} of
$u$. Then, $Du$ belongs to the class $M^{b}(\Omega,\R^{d})$ of $\R^{d}$-valued bounded
Radon measure on $\Omega$, and we either write $\abs{Du}(\Sigma)$ or
$\int_\Sigma \vert Du \vert$ to denote the \emph{total variation
  measure} of $Du$. The space $BV(\Sigma)$ equipped with the norm
\begin{math}
  \norm{u}_{BV(\Sigma)}:= \norm{u}_{L^1(\Sigma)} +
  \abs{Du}(\Sigma)
\end{math}
forms a Banach space. Further, let
\begin{displaymath}
  X_{1}(\Sigma)=\Big\{z\in L^{\infty}(\Sigma,\R^{d})\,\Big\vert\Big.\;
  \textrm{div}(z):=\sum_{i=1}^{d}D_{i}z\in L^{1}(\Sigma)\Big\},
\end{displaymath}
$\textrm{sign}_0(s)$, ($s\in \R$), is the classical sign function with the additional
property that $\textrm{sign}_0(0)=0$, and for every $k>0$, $T_k(s):=
\left[k - [k - \abs{s}\right]^+ \textrm{sign}_0(s)$, ($s\in \R$).\medskip

{\bfseries The Neumann total variational flow operator.} In
\cite{MR1799898} (see also \cite{MR2033382}), the negative \emph{total
  variational flow operator} (\emph{$1$-Laplace operator})
$-\Delta_{1}^{\! N}$ in $L^{1}(\Sigma)$ equipped with Neumann boundary
conditions was introduced by
\begin{displaymath}
  -\Delta^{\! N}_{1}=\Bigg\{(u,v)\in L^{1}(\Sigma)\times
  L^{1}(\Sigma)\,\Bigg\vert\,
  \begin{array}[c]{l}
T_{k}(u)\in BV(\Sigma)\,\forall\,k>0\text{ \& }
\exists\;z\in X_{1}(\Sigma)\\
 \text{such that $\norm{z}_{\infty}\le 1$ \&~\eqref{eq:NTF} holds}
  \end{array}
\Bigg\},
\end{displaymath}
where
\begin{equation}
  \label{eq:NTF}
  \begin{cases}
     v=-\textrm{div}(z)\qquad\text{in $\mathcal{D}'(\Sigma)$, and}&\\[7pt]
    \displaystyle\int_{\Sigma}\left(\xi-T_{k}(u\right) v\,\dx\le
    \int_{\Sigma}z\cdot D\xi\,\dx-\int_{\Sigma}\abs{DT_{k}(u)}&
  \end{cases}
\end{equation}
for every $\xi\in W^{1,1}(\Sigma)\cap L^{\infty}$ and all
$k>0$. Moreover, the negative Neumann $1$-Laplace
operator $-\Delta^{\! N}_{1}$ is $m$-completely accretive in $L^{1}(\Sigma)$
with dense domain. Therefore, under the
hypotheses~\eqref{property:one-Charatheodory}--\eqref{eq:2}, the
operator $-\Delta^{\! N}_{1}+F$ is $\omega$-quasi $m$-completely
accretive on $L^{1}(\Sigma)$ (cf~\cite{CoulHau2017}). Now, it is not difficult to see that
$-\Delta^{\! N}_{1}$ is homogeneous of order zero and $0\in -\Delta^{\! N}_{1}0$. Thus, by
Theorem~\ref{thm:regularity-of-complete-acc-homogen-operators} and
Theorem~\ref{thm:Lipschitz}, we can state the following regularity result.

\begin{corollary}\label{cor:TVN}
  For every $1\le q<\infty$ and $u_{0}\in L^{\psi}(\Sigma)$
  (respectively $u_{0}\in L^{1\cap\infty}(\Sigma)$ if $q=\infty$), the unique mild solution $u$
  of problem
  \begin{equation}
    \label{eq:NTF-Cauchy}
    \begin{cases}
     \displaystyle \frac{\td u}{\dt} -\textrm{div}\left(\frac{D u}{\abs{D u}}\right)+ f(x,u)=0
      & \text{on $\Sigma\times (0,+\infty)$,}\\[7pt]
      \hspace{3.95cm}D_{\nu}u =0 & \text{on $\partial\Sigma\times (0,+\infty)$,}\\[7pt]
      \hspace{3.9cm}u(0)=u_{0} & \text{on $\Sigma\times \{t=0\}$,}\\
    \end{cases}
  \end{equation}
  is a strong solution satisfying~\eqref{eq:6-F-Lq}.
  Moreover, if $f \equiv 0$, then either for every
  $1\le \psi\equiv p\le\infty$ or $N$-function $\psi$ and every
  $u_{0}\in L^{\psi}(\Sigma)$, the unique mild solution $u$ of
  problem~\eqref{eq:NTF-Cauchy} satisfies~\eqref{eq:6} and~\eqref{eq:4}.
\end{corollary}

\begin{remark}[{\bfseries The Dirichlet boundary case}]
  In~\cite{MR1814993} (cf~\cite{MR2033382}), existence and uniqueness
  of the the parabolic initial boundary-value problem
 \begin{equation}\label{eq:DirichletTV}
\begin{cases}
     \displaystyle \frac{\td u}{\dt} -\textrm{div}\left(\frac{D u}{\abs{D u}}\right)=0
      & \text{on $\Sigma\times (0,+\infty)$,}\\[7pt]
      \hspace{2.9cm}u =\varphi & \text{on $\partial\Sigma\times (0,+\infty)$,}\\[7pt]
      \hspace{2.4cm}u(0)=u_{0} & \text{on $\Sigma\times \{t=0\}$,}\\
    \end{cases}
 \end{equation}
 associated with the total variational flow equipped with
 (inhomogeneous) Diri\-chlet boundary conditions was established. For
 every boundary term $\varphi\in L^{1}(\Sigma)$, the negative \emph{Dirichlet
   total variational flow operator ($1$-Laplace operator)}
 $\Delta_{1}^{\! D}u:=\textrm{div}\left(\frac{D u}{\abs{D
       u}}\right)$ is $m$-completely accretive in $L^{1}(\Sigma)$. But only in the homogeneous case $\varphi\equiv
 0$, the operator $\Delta^{\! D}_{1}$ is homogeneous of order zero. Thus, the same
 statement as given in Corollary~\ref{cor:TVN} holds in the Dirichlet
 case with $\varphi \equiv 0$.
\end{remark}

%
%
\subsection{Decay estimates of the nonlocal total variational Flow}
In this very last section, we consider for $0<s<1$, the perturbed fractional
$1$-Laplace operator
\begin{displaymath}
  Au:=\textrm{PV}\int_{\Sigma}\frac{(u(y)-u(x))}{\abs{u(y)-u(x)}}\, \frac{\dy}{\abs{x-y}^{d+s}}+f(x,u)
\end{displaymath}
equipped with either Dirichlet on a domain $\Sigma$ in $\R^{d}$ or
with or vanishing conditions if $\Sigma=\R^{d}$, $d\ge 1$.\medskip

For $0<s<1$, let $\W_{0}^{s,1}(\Sigma)$ be the Banach space given by
\begin{displaymath}
  \W_{0}^{s,1}(\Sigma)=\Big\{u\in L^{1}(\Sigma)\,\Big\vert\,
          [u]_{s,1}<\infty\text{ and }u=0\text{ a.e. on }\R^{d}\setminus\Sigma\Big\}
\end{displaymath}
equipped with the norm $\norm{\cdot}_{\W_{0}^{s,1}}:=\norm{\cdot}_{1}+[\cdot]_{s,1}$, where
\begin{displaymath}
  [u]_{s,1}:=\int_{\R^{d}}\int_{\R^{d}}\frac{\abs{u(x)-u(y)}}{\abs{x-y}^{d+s}}\dy\dx\qquad\text{for
  every $u\in \W_{0}^{s,1}(\Sigma)$.}
\end{displaymath}
Further, let $B_{L^{\infty}_{as}}$ denote the closed unit ball of all
\emph{anti-symmetric} $\eta \in L^{\infty}(\R^{d}\times\R^{d})$, that is,
\begin{displaymath}
  \eta(x,y)=-\eta(y,x)\qquad\text{ for a.e. }(x,y)\in
  \R^{d}\times\R^{d}\quad\text{ and }\norm{\eta}_{\infty}\le 1.
\end{displaymath}
Then, it was shown in~\cite[Section~3.]{MR3491533} that the
\emph{fractional Dirichlet $1$-Laplace operator
  $(-\Delta^{D}_{1})^{s}$} in $L^{2}(\Sigma)$ can be realized by (the graph)
\begin{displaymath}
  (-\Delta^{D}_{1})^{s}=\Big\{(u,v)\in L^{2}(\Sigma)\times
  L^{2}(\Sigma)\,\Big\vert\,
  u\in \W_{0}^{s,1}(\Sigma)\text{ \& }\exists\;\eta\in
  B_{L^{\infty}_{as}}\text{ s.t.~\eqref{eq:34} holds}\Big\}
\end{displaymath}
where
\begin{equation}
  \label{eq:34}
  \begin{cases}
     \eta(x,y)\in \textrm{sign}(u(x)-u(y))\qquad\text{for
      a.e. $(x,y)\in \R^{d}\times\R^{d}$ and}&\\[7pt]
    \displaystyle \frac{1}{2}\int_{\R^{d}}\int_{\R^{d}}\frac{\eta(x,y)
      (\xi(x)-\xi(y))}{\abs{y-x}^{d+s}}\;\dy\dx
    =\int_{\Sigma}v(x)\,\xi(x)\,\dx&\\[7pt]
    \hspace{5.6cm}\text{for all $\xi\in \W_{0}^{s,1}(\Sigma)\cap L^{2}(\Sigma)$,}&
  \end{cases}
\end{equation}
and $(-\Delta^{D}_{1})^{s}$ is $m$-completely accretive in
$L^{2}(\Sigma)$ with dense domain $D((-\Delta^{D}_{1})^{s})$ in
$L^{2}(\Sigma)$. One immediately sees that $(-\Delta^{D}_{1})^{s}$ is
homogeneous of order zero and $0\in (-\Delta^{D}_{1})^{s}0$. Moreover,
under the hypotheses~\eqref{property:one-Charatheodory}--\eqref{eq:2},
the operator $(-\Delta^{D}_{1})^{s}+F$ is $\omega$-quasi
$m$-completely accretive on $L^{2}(\Sigma)$. Thus, by
Theorem~\ref{thm:Lipschitz}, we have the following regularity result.


\begin{corollary}\label{cor:nonlocal-TVF}
  For every $1\le q<\infty$ and $u_{0}\in L^{\psi}(\Sigma)$
  (respectively $u_{0}\in L^{1\cap\infty}(\Sigma)$ if $q=\infty$), the unique mild solution $u$
  of problem
  \begin{equation}
    \label{eq:35}
    \begin{cases}
     \displaystyle \frac{\td
       u}{\dt}+\textrm{PV}\int_{\Sigma}\frac{u(x)-u(y)}{\abs{u(y)-u(x)}}
     \frac{\dy}{\abs{x-y}^{d+s}}+ f(x,u)=0
      & \text{on $\Sigma\times (0,+\infty)$,}\\[7pt]
      \hspace{7.3cm}u=0 & \text{on $\partial\Sigma\times (0,+\infty)$,}\\[7pt]
      \hspace{6.8cm}u(0)=u_{0} & \text{on $\Sigma\times \{t=0\}$,}\\
    \end{cases}
  \end{equation}
  is a strong solution satisfying~\eqref{eq:6-F-Lq}. Moreover, if
  $f \equiv 0$, then either for every $1\le \psi\equiv p\le\infty$ or
  $N$-function $\psi$ and every $u_{0}\in L^{\psi}(\Sigma)$, the
  unique mild solution $u$ of problem~\eqref{eq:NTF-Cauchy}
  satisfies~\eqref{eq:6} and~\eqref{eq:4}.
\end{corollary}

Further, since for every $0<s<1$ and $d\ge 1$, the following
\emph{(fractional) Sobolev inequality} (cf~\cite[Theorem~14.29]{MR2527916})
\begin{displaymath}
  \norm{u}_{\frac{d}{d-s}}\le C\,[u]_{s,1}\qquad\text{for every $u\in \W_{0}^{s,1}(\Sigma)$,}
\end{displaymath}
holds for a constant $C=C(d,s)>0$, we have that the fractional
Dirichlet $1$-Laplace operator $A=(-\Delta^{D}_{1})^{s}$ satisfies the
following abstract Sobolev inequality
\begin{displaymath}
  \norm{u}_{r}^{\sigma}\le C\,[u,v]_{2}\qquad\text{for every $(u,v)\in
    A$}
\end{displaymath}
with parameters $r=\frac{d}{d-s}>1$ and $\sigma=1$, where
$[\cdot,\cdot]_{2}$ denote the $L^{2}$-inner product. Thus,
by~\cite[Theorem~1.2]{CoulHau2017}, the semigroup $\{T_{t}\}_{t\ge 0}$
generated by $-((-\Delta^{D}_{1})^{s}+F)$ satisfies the following
\emph{$L^{2}$-$L^{\frac{d}{d-s}}$-regularity estimate}
\begin{displaymath}
  \norm{T_{t}u}_{\frac{d}{d-s}}\le \frac{C}{2}\,t^{-1}\,e^{3\omega
    t}\norm{u_{0}}_{2}^{2}\qquad\text{for every $t>0$}
\end{displaymath}
and $u_{0}\in L^{2}(\Sigma)$. Furthermore, for every
$q> \frac{d-s}{s}$, one has
\begin{displaymath}
  \norm{T_{t}u}_{\infty}\le \tilde{C}
  \,t^{-\frac{d-s}{s(q+1)-d}}\,e^{\left(\frac{d-s}{s(q+1)-d}+1\right) \omega t}\,\norm{u_{0}}_{\frac{dq}{d-s}}^{\frac{d-s}{s(q+1)-d}\frac{sq}{d-s}}
  \qquad\text{for every $t>0$}
\end{displaymath}
and $u_{0}\in L^{\frac{dq}{d-s}}(\Sigma)$. Thus, by
Corollary~\ref{cor:1} and inequality~\eqref{eq:6-F-Lq} for $q=\infty$, we
also have the following estimate.

 \begin{corollary}\label{cor:Decay-nonlocal-DTF}
  Let $d\ge 1$, $0<s<1$ and $q>\frac{d-s}{s}$. Then
  for every $u_{0}\in L^{\frac{dq}{d-s}}(\Sigma)$, the unique solution $u$ of~\eqref{eq:35} satisfies
  \begin{displaymath}
    \lnorm{\frac{\td u}{\dt}}_{\infty}\le  \tilde{C}\,C_{\omega}(t/2)
    2^{\frac{d-s}{s(q+1)-d}+2} e^{\omega \left(\frac{d-s}{2s(q+1)-d}+1\right)
      t}\,\frac{\norm{u_{0}}_{\psi}^{\frac{d-s}{s(q+1)-d}\frac{sq}{d-s}}}{t^{\frac{d-s}{s(q+1)-d}+1}}
    \qquad\text{for every $t>0$,}
  \end{displaymath}
  where $C_{\omega}(t)$ is the constant in~\eqref{eq:6-F-Lq}.
\end{corollary}

%
%


\end{document}